\documentclass{siamonline1116}

\usepackage{lipsum}
\usepackage{amsfonts}
\usepackage{graphicx}
\usepackage{epstopdf}
\usepackage{algorithmic}
\usepackage{amsmath}
\usepackage{amssymb}
\usepackage{latexsym}
\usepackage[caption=false]{subfig}
\ifpdf
  \DeclareGraphicsExtensions{.eps,.pdf,.png,.jpg}
\else
  \DeclareGraphicsExtensions{.eps}
\fi

\numberwithin{theorem}{section}

\newcommand{\TheTitle}{Linearly-Recurrent Autoencoder Networks for Learning Dynamics}
\newcommand{\TheAuthors}{S. E. Otto and C. W. Rowley}

\headers{\TheTitle}{\TheAuthors}

\title{{\TheTitle}\thanks{Submitted to the editors \today.
\funding{This work was funded by ARO award W911NF-17-0512 and DARPA.}}}

\author{
  Samuel E. Otto\thanks{Mechanical and Aerospace Engineering, Princeton University, Princeton, NJ
    (\email{sotto@princeton.com}).}
  \and
  Clarence W. Rowley\thanks{Mechanical and Aerospace Engineering, Princeton University, Princeton, NJ (\email{cwrowley@princeton.edu}).}
}

\usepackage{amsopn}

\DeclareMathOperator*{\minimize}{\min\!imize\enskip}
\DeclareMathOperator{\Tr}{Tr}
\DeclareMathOperator{\diag}{diag}
\DeclareMathOperator{\vspan}{span}
\DeclareMathOperator{\elu}{elu}
\DeclareMathOperator*{\Expectation}{\!\mathbb{E}\enskip}

\newsiamremark{example}{Example}

\ifpdf
\hypersetup{
  pdftitle={\TheTitle},
  pdfauthor={\TheAuthors}
}
\fi

\usepackage{xr}
\begin{document}

\maketitle

% REQUIRED
\begin{abstract}
  This paper describes a method for learning low-dimensional approximations of
  nonlinear dynamical systems, based on neural-network approximations of the
  underlying Koopman operator. Extended Dynamic Mode Decomposition (EDMD)
  provides a useful data-driven approximation of the Koopman operator for
  analyzing dynamical systems. This paper addresses a fundamental problem
  associated with EDMD: a trade-off between representational capacity of the
  dictionary and over-fitting due to insufficient data. A new neural network
  architecture combining an autoencoder with linear recurrent dynamics in the
  encoded state is used to learn a low-dimensional and highly informative
  Koopman-invariant subspace of observables. A method is also presented for
  balanced model reduction of over-specified EDMD systems in feature
  space. Nonlinear reconstruction using partially linear multi-kernel regression
  aims to improve reconstruction accuracy from the low-dimensional state when
  the data has complex but intrinsically low-dimensional structure. The
  techniques demonstrate the ability to identify Koopman eigenfunctions of the
  unforced Duffing equation, create accurate low-dimensional models of an
  unstable cylinder wake flow, and make short-time predictions of the chaotic
  Kuramoto-Sivashinsky equation.
\end{abstract}

% REQUIRED
\begin{keywords}
  nonlinear systems, high-dimensional systems, reduced-order modeling, neural
  networks, data-driven analysis, Koopman operator
\end{keywords}

% REQUIRED
\begin{AMS}
  65P99, 37M10, 37M25, 47-04, 47B33
\end{AMS}

\section{Introduction}
The Koopman operator first introduced in \cite{Koopman1931} describes how
Hilbert space functions on the state of a dynamical system evolve in time. These
functions, referred to as observables, may correspond to measurements taken
during an experiment or the output of a simulation. This makes the Koopman
operator a natural object to consider for data-driven analysis of dynamical
systems. Such an approach is also appealing because the Koopman operator is
linear, though infinite dimensional, enabling the concepts of modal analysis for
linear systems to be extended to dynamics of observables in nonlinear
systems. Hence, the invariant subspaces and eigenfunctions of the Koopman
operator are of particular interest and provide useful features for describing
the system if they can be found. For example, level sets of Koopman
eigenfunctions may be used to form partitions of the phase space into ergodic
sets along with periodic and wandering chains of sets \cite{Budisic2012}. They
allow us to parameterize limit cycles and tori as well as their basins of
attraction. The eigenvalues allow us to determine the stability of these
structures and the frequencies of periodic and quasiperiodic attractors
\cite{Mezic2017}. Furthermore, by projecting the full state as an observable
onto the eigenfunctions of the Koopman operator, it is decomposed into a linear
superposition of components called Koopman modes which each have a fixed
frequency and rate of decay. Koopman modes therefore provide useful coherent
structures for studying the system's evolution and dominant pattern-forming
behaviors. This has made the Koopman operator a particularly useful object of
study for high-dimensional spatiotemporal systems like unsteady fluid dynamics
beginning with the work of Mezi\'{c} on spectral properties of dynamical systems
\cite{Mezic2005} then Rowley \cite{ROWLEY2009} and Schmid \cite{SCHMID2010} on
the Dynamic Mode Decomposition (DMD) algorithm. Rowley, recognizing that DMD
furnishes an approximation of the Koopman operator and its modes, applied the
technique to data collected by simulating a jet in a crossflow. The Koopman
modes identified salient patterns of spatially coherent structure in the flow
which evolved at fixed frequencies.

The Extended Dynamic Mode Decomposition (EDMD) \cite{Williams2015} is an
algorithm for approximating the Koopman operator on a dictionary of observable
functions using data. If a Koopman-invariant subspace is contained in the span
of the observables included in the dictionary, then as long as enough data is
used, the representation on
this subspace will be exact. EDMD is a Galerkin method with a particular data-driven inner product as long as enough data is used. 
Specifically, this will be true as long as the rank of the data matrix is the same as the dimension of the subspace spanned by the (nonlinear) observables \cite{Rowley2017_Kavli}. 
However, the choice of dictionary is ad hoc, and it is often not clear how to
choose a dictionary that is sufficiently rich to
span a useful Koopman-invariant subspace. One might then be tempted to consider
a very large dictionary, with enough capacity to represent any complex-valued function on the
state space to within an $\epsilon$ tolerance. 
However, such a dictionary has
combinatorial growth with the dimension of the state space and would be enormous
for even modestly high-dimensional problems.

One approach to mitigate the cost
of large or even infinite dictionaries is to formulate EDMD as a kernel method
referred to as KDMD \cite{Kevrekidis2016}. However, we are still essentially
left with the same problem of deciding which kernel function to
use. Furthermore, if the kernel or EDMD feature space is endowed with too much
representational capacity (a large dictionary), the algorithm will over-fit the
data (as we shall demonstrate with a toy problem in
\Cref{ex:EDMD_counterexample}).
EDMD and KDMD also identify a number of eigenvalues, eigenfunctions, and
modes which grows with the size of the dictionary. If we want to build reduced
order models of the dynamics, a small collection of salient modes or a
low-dimensional Koopman invariant subspace must be identified. It is worth
mentioning two related algorithms for identifying low-rank approximations of the
Koopman operator. Optimal Mode Decomposition (OMD) \cite{Wynn2013} finds the
optimal orthogonal projection subspace of user-specified rank for approximating
the Koopman operator. Sparsity-promoting DMD \cite{Jovanovic2014} is a
post-processing method which identifies the optimal amplitudes of Koopman modes
for reconstructing the snapshot sequence with an $\ell^1$ penalty. The
sparsity-promoting penalty picks only the most salient Koopman modes to have
nonzero amplitudes. Another related scheme is Sparse Identification of Nonlinear
Dynamics (SINDy) \cite{Brunton2016} which employs a sparse regression penalty on
the number of observables used to approximate nonlinear evolution equations. By
forcing the dictionary to be sparse, the over-fitting problem is reduced.

In this paper, we present a new technique for learning a very small collection
of informative observable functions spanning a Koopman invariant subspace
from data. Two neural networks in an architecture similar to an under-complete
autoencoder \cite{Goodfellow2016} represent the collection of observables
together with a nonlinear reconstruction of the full state from these
features. A learned linear transformation evolves the function values in time as
in a recurrent neural network, furnishing our approximation of the Koopman
operator on the subspace of observables. This approach differs from recent
efforts that use neural networks to learn dictionaries for EDMD \cite{Yeung2017,
  Li2017} in that we employ a second neural network to reconstruct the full
state. Ours and concurrent approaches utilizing nonlinear decoder neural networks \cite{Takeishi2017, Lusch2017} enable learning of very small sets of features that carry rich information about the state and evolve linearly in time. Previous methods for data-driven analysis based
on the Koopman operator utilize linear state reconstruction via the Koopman
modes. Therefore they rely on an assumption that the full state observable is in
the learned Koopman invariant subspace. 
Nonlinear reconstruction is advantageous since it relaxes this strong
assumption, allowing recent techniques to recover more information about the state from fewer
observables. By minimizing the
state reconstruction error over several time steps into the future, our architecture
aims to detect highly observable features even if they have small
amplitudes. This is the case in non-normal linear systems, for instance as arise
in many fluid flows (in particular, shear flows \cite{SchmidHenn}), in which small disturbances can siphon energy from mean flow gradients and
excite larger-amplitude modes. The underlying philosophy of our approach is similar to Visual Interaction Networks (VINs) \cite{Watters2017} that learn physics-based dynamics models for encoded latent variables.

Deep neural networks have gained attention over the last decade due to their
ability to efficiently represent complicated functions learned from data. Since each
layer of the network performs simple operations on the output of the previous
layer, a deep network can learn and represent functions
corresponding to high-level or abstract features. For example, your visual
cortex assembles progressively more complex information sequentially from
retinal intensity values to edges, to shapes, all the way up to the facial
features that let you recognize your friend. By contrast, shallow networks ---
though still universal approximators --- require exponentially more parameters
to represent classes of natural functions like polynomials \cite{Rolnick2017,
  Lin2016} or the presence of eyes in a photograph. Function approximation using
linear combinations of preselected dictionary elements is somewhat analogous to
a shallow neural network where capacity is built by adding more functions. We
therefore expect deep neural networks to represent certain complex nonlinear
observables more efficiently than a large, shallow dictionary. Even with the
high representational capacity of our deep neural networks, the proposed
technique is regularized by the small number of observables we learn and is
therefore unlikely to over-fit the data.

We also present a technique for constructing reduced order models in nonlinear
feature space from over-specified KDMD models. Recognizing that the systems
identified from data by EDMD/KDMD can be viewed as state-space systems where the
output is a reconstruction of the full state using Koopman modes, we use
Balanced Proper Orthogonal Decomposition (BPOD) \cite{Rowley2005} to construct a
balanced reduced-order model. The resulting model consists of only those
nonlinear features that are most excited and observable over a finite time
horizon. Nonlinear reconstruction of the full state is introduced in order to
account for complicated, but intrinsically low-dimensional data. In this way,
the method is analogous to an autoencoder where the nonlinear decoder is learned
separately from the encoder and dynamics.

Finally, the two techniques we introduce are tested on a range of example
problems. We first investigate the eigenfunctions learned by the autoencoder and
the KDMD reduced-order model by identifying and parameterizing basins of
attraction for the unforced Duffing equation. The prediction accuracy of the
models is then tested on a high-dimensional cylinder wake flow problem. Finally,
we see if the methods can be used to construct reduced order models for the
short-time dynamics of the chaotic Kuramoto-Sivashinsky equation. Several
avenues for future work and extensions of our proposed methods are discussed in
the conclusion.

\section{Extended Dynamic Mode Decomposition}
Before discussing the new method for approximating the Koopman operator, it will
be beneficial to review the formulation of Extended Dynamic Mode Decomposition
(EDMD) \cite{Williams2015} and its kernel variant KDMD
\cite{Kevrekidis2016}. Besides providing the context for developing the new
technique, it will be useful to compare our results to those obtained using reduced order
KDMD models.

\subsection{The Koopman operator and its modes}
Consider a discrete-time autonomous dynamical system on the state space
$\mathcal{M} \subset \mathbb{R}^n$ given by the function
$\mathbf{x}_{t+1} = \mathbf{f}(\mathbf{x}_t)$. Let $\mathcal{F}$ be a Hilbert
space of complex-valued functions on~$\mathcal{M}$.  We refer to elements of $\mathcal{F}$
as observables.  The Koopman
operator acts on an observable $\psi\in \mathcal{F}$ by composition with the dynamics:
\begin{equation} \label{eqn:Koopman_def}
  \mathcal{K}\psi = \psi \circ \mathbf{f}.
\end{equation}
It is easy to see that the Koopman operator is linear; however, the Hilbert
space $\mathcal{F}$ on which it acts is often infinite dimensional.%
\footnote{One must also be careful about the choice of the space $\mathcal{F}$,
  since $\psi\circ\mathbf{f}$ must also lie in $\mathcal{F}$ for any
  $\psi\in\mathcal{F}$.  It is common, especially in the ergodic theory
  literature, to assume that $\mathcal{M}$ is a measure space and $\mathbf{f}$ is measure
  preserving.  In this case, this difficulty goes away: one lets
  $\mathcal{F}=L^2(\mathcal{M})$, and since $\mathbf{f}$ is measure preserving, it follows
  that $\mathcal{K}$ is an isometry.}
Since the operator~$\mathcal{K}$ is linear, it may have eigenvalues and
eigenfunctions.
If a given observable lies within the span of these eigenfunctions, then we can
predict the time evolution of the observable's values, as the state evolves according to the dynamics. Let
$\mathbf{g}:\mathcal{M}\rightarrow\mathbb{C}^{N_0}$ be a vector-valued
observable whose components are in the span of the Koopman eigenfunctions. The
vector-valued coefficients needed to reconstruct $\mathbf{g}$ in a Koopman
eigenfunction basis are called the Koopman modes associated with $\mathbf{g}$.

In particular, the dynamics of the original system can be recovered by taking
the observable~$\mathbf{g}$ to be the full-state observable defined by
$\mathbf{g}(\mathbf{x}) = \mathbf{x}$.  Assume $\mathcal{K}$ has eigenfunctions
$\{\varphi_1,\ldots,\varphi_K\}$ with
corresponding eigenvalues $\{\mu_1,\ldots,\mu_K\}$, and suppose the components
of the vector-valued function~$\mathbf{g}$ lie within the span of
$\{\varphi_k\}$.  The Koopman modes $\boldsymbol{\xi}_k$ are then defined by
\begin{equation}
  \label{eqn:Koopman_modes}
  \mathbf{x} = \sum_{k=1}^K \boldsymbol{\xi}_k \varphi_k(\mathbf{x}),
\end{equation}
from which we can recover the evolution of the state, according to
\begin{equation}
  \mathbf{f}^t(\mathbf{x}) = \sum_{k=1}^K \boldsymbol{\xi}_k \mu_k^t \varphi_k(\mathbf{x}).
\end{equation}
The entire orbit of an
initial point $\mathbf{x}_0$ may thus be determined by evaluating the eigenfunctions at
$\mathbf{x}_0$ and evolving the coefficients $\boldsymbol{\xi}_k$ in time by
multiplying by the eigenvalues. The eigenfunctions $\varphi_k$ are intrinsic features of
the dynamical system which decompose the state dynamics into a linear
superposition of autonomous first-order systems.
The Koopman modes $\boldsymbol{\xi}_k$ depend on the coordinates we use to
represent the dynamics, and allow us to reconstruct the dynamics in those coordinates.

\subsection{Approximating Koopman on an explicit dictionary with EDMD}
The aim of EDMD is to approximate the Koopman operator using data snapshot pairs
taken from the system
$\left\lbrace \left(\mathbf{x}_j, \mathbf{y}_j \right)\right\rbrace_{j=1}^M$
where $\mathbf{y}_j = \mathbf{f}(\mathbf{x}_j)$. For convenience, we organize
these data into matrices
\begin{equation}
  \mathbf{X} =
  \begin{bmatrix}
    \mathbf{x}_1 & \mathbf{x}_2 & \cdots & \mathbf{x}_M
  \end{bmatrix}, \qquad
  \mathbf{Y} =
  \begin{bmatrix}
    \mathbf{y}_1 & \mathbf{y}_2 & \cdots & \mathbf{y}_M
  \end{bmatrix}.
\end{equation}
Consider a finite dictionary of observable functions
$\mathcal{D} = \left\lbrace \psi_i:\mathcal{M} \rightarrow \mathbb{C}
\right\rbrace _{i=1}^N$ that span a subspace
$\mathcal{F}_{\mathcal{D}} \subset \mathcal{F}$. EDMD approximates the Koopman
operator on $\mathcal{F}_{\mathcal{D}}$ by minimizing an empirical error when
the Koopman operator acts on elements $\psi \in
\mathcal{F}_{\mathcal{D}}$. Introducing the feature map
\begin{equation}
  \mathbf{\Psi}(\mathbf{x}) =
  \begin{bmatrix}
    \psi_1(\mathbf{x}) & \psi_2(\mathbf{x}) & \cdots & \psi_N(\mathbf{x})
  \end{bmatrix}^*,
\end{equation}
where $(\cdot)^*$ is the complex conjugate transpose, we may
succinctly express elements in the dictionary's span as a linear combination with
coefficients $\mathbf{a}$:
\begin{equation}
  \label{eqn:observable_from_features}
  \psi_{\mathbf{a}} = \mathbf{\Psi}^* \mathbf{a}.
\end{equation}
EDMD represents an approximation of the Koopman operator as a matrix
$\mathbf{K}:\mathbb{C}^N \rightarrow \mathbb{C}^N$ that updates the coefficients
in the linear combination \cref{eqn:observable_from_features} to approximate the
new observable $\mathcal{K} \psi_{\mathbf{a}}$ in the span of the dictionary. Of
course we cannot expect the span of our dictionary to be an invariant subspace,
so the approximation satisfies
\begin{equation}
  \label{eqn:Koopman_approx_on_observable}
  \mathcal{K} \psi_{\mathbf{a}} = \mathbf{\Psi}^* \mathbf{K} \mathbf{a} + r,
\end{equation}
where $r\in\mathcal{F}$ is a residual that we wish to minimize in some sense, by appropriate
choice of the matrix~$\mathbf{K}$.
The values of the Koopman-updated observables are known at each of the
data points
$\mathcal{K} \psi_{\mathbf{a}}(\mathbf{x}_j) = \psi_{\mathbf{a}}\circ
\mathbf{f}(\mathbf{x}_j) = \psi_{\mathbf{a}}(\mathbf{y}_j)$, allowing us to
define an empirical error of the approximation in terms of the
residuals at these data points. Minimizing this error yields the EDMD matrix $\mathbf{K}$. The empirical error on a single
observable in $\mathcal{F}_{\mathcal{D}}$ is given by
\begin{equation}
  \label{eqn:EDMD_single_observable_error}
  \begin{split}
    J(\psi_{\mathbf{a}}) =& \sum_{i=1}^M \left\vert \psi_{\mathbf{a}}(\mathbf{y}_i) - \mathbf{\Psi}(\mathbf{x}_i)^* \mathbf{K} \mathbf{a} \right\vert^2 \\
    =& \sum_{i=1}^M \left\vert \left( \mathbf{\Psi}(\mathbf{y}_i)^* - \mathbf{\Psi}(\mathbf{x}_i)^* \mathbf{K} \right) \mathbf{a} \right\vert^2
  \end{split}
\end{equation}
and the total empirical error on a set of observables
$\left\lbrace \psi_{\mathbf{a}_j} \right\rbrace_{j=1}^{N'}$, $N' \geq N$
spanning $\mathcal{F}_{\mathcal{D}}$ is given by
\begin{equation}
  \label{eqn:EDMD_total_error}
  J = \sum_{j=1}^{N'} \sum_{i=1}^M \left\vert \left( \mathbf{\Psi}(\mathbf{y}_i)^* - \mathbf{\Psi}(\mathbf{x}_i)^* \mathbf{K} \right) \mathbf{a}_j \right\vert^2.
\end{equation}
Regardless of how the above observables are chosen, the matrix~$\mathbf{K}$ that
minimizes \cref{eqn:EDMD_total_error} is given by
\begin{equation}
  \label{eqn:EDMD_solution}
  \mathbf{K} = \mathbf{G}^{+}\mathbf{A}, \qquad
  \mathbf{G} = \frac{1}{M} \sum_{i=1}^M \mathbf{\Psi}(\mathbf{x}_i)\mathbf{\Psi}(\mathbf{x}_i)^*, \qquad
  \mathbf{A} = \frac{1}{M} \sum_{i=1}^M \mathbf{\Psi}(\mathbf{x}_i)\mathbf{\Psi}(\mathbf{y}_i)^*,
\end{equation}
where $(\cdot)^{+}$ denotes the Moore-Penrose pseudoinverse of a matrix.

The EDMD solution \cref{eqn:EDMD_solution} requires us to evaluate the entries
of $\mathbf{G}$ and $\mathbf{A}$ and compute the pseudoinverse of
$\mathbf{G}$. Both matrices have size $N\times N$ where $N$ is the number of
observables in our dictionary. In problems where the state dimension is large,
as it is in many fluids datasets coming from experimental or simulated flow
fields, a very large number of observables is needed to achieve even modest
resolution on the phase space. The problem of evaluating and storing the
matrices needed for EDMD becomes intractable as $N$ grows large. However, the
rank of these matrices does not exceed $\min \left\lbrace M,N
\right\rbrace$. The kernel DMD method provides a way to compute an EDMD-like
approximation of the Koopman operator using kernel matrices whose size scales
with the number of snapshot pairs $M^2$ instead of the number of features
$N^2$. This makes it advantageous for problems where the state dimension is
greater than the number of snapshots or where very high resolution of the
Koopman operator on a large dictionary is needed.

\subsection{Approximating Koopman on an implicit dictionary with KDMD}
KDMD can be derived by considering the data matrices %
\begin{equation}
  \label{eqn:feature_matrices}
  \mathbf{\Psi}_{\mathbf{X}} =
  \begin{bmatrix}
    \mathbf{\Psi}(\mathbf{x}_1) & \mathbf{\Psi}(\mathbf{x}_2) & \cdots & \mathbf{\Psi}(\mathbf{x}_M)
  \end{bmatrix}, \qquad
  \mathbf{\Psi}_{\mathbf{Y}} =
  \begin{bmatrix}
    \mathbf{\Psi}(\mathbf{y}_1) & \mathbf{\Psi}(\mathbf{y}_2) & \cdots & \mathbf{\Psi}(\mathbf{y}_M)
  \end{bmatrix},
\end{equation}
in feature space i.e., after applying the now only hypothetical feature map $\mathbf{\Psi}$ to the snapshots. 
We will see that the final results of this approach make reference only to inner products
$\mathbf{\Psi}(\mathbf{x})^*\mathbf{\Psi}(\mathbf{z})$ which will be defined
using a suitable non-negative definite kernel function
$k(\mathbf{x},\mathbf{z})$. Choice of such a kernel function implicitly defines
the corresponding dictionary via Mercer's theorem. By employing simply-defined
kernel functions, the inner products are evaluated at a lower computational cost
than would be required to evaluate a high or infinite dimensional feature map
and compute inner products in the feature space explicitly. 

The total empirical error for EDMD \cref{eqn:EDMD_total_error} can be written as the Frobenius norm
\begin{equation}
  \label{eqn:EDMD_Frobenius_error}
  J = \left\Vert \left(\mathbf{\Psi}_{\mathbf{Y}}^* - \mathbf{\Psi}_{\mathbf{X}}^* \mathbf{K} \right) \mathbf{A} \right\Vert_F^2, \quad \mbox{where} \quad
  \mathbf{A} =
  \begin{bmatrix}
    \mathbf{a}_1 & \mathbf{a}_2 & \cdots & \mathbf{a}_{N'}
  \end{bmatrix}.
\end{equation}
Let us consider an economy sized SVD
$\mathbf{\Psi}_{\mathbf{X}} = \mathbf{U}\mathbf{\Sigma}\mathbf{V}^*$, the
existence of which is guaranteed by the finite rank $r$ of our feature data
matrix. In \cref{eqn:EDMD_Frobenius_error} we see that any components of the
range $\mathcal{R}(\mathbf{K})$ orthogonal to
$\mathcal{R}(\mathbf{\Psi}_{\mathbf{X}})$ are annihilated by
$\mathbf{\Psi}_{\mathbf{X}}^*$ and cannot be inferred from the data. We
therefore restrict the dictionary to those features which can be represented in
the range of the feature space data
$\mathcal{F}_{\mathcal{D}} = \left\lbrace \psi_{\mathbf{a}} =
  \mathbf{\Psi}^*\mathbf{a} : \mathbf{a}\in\mathcal{R}(\mathbf{U})
\right\rbrace$ and represent
$\mathbf{K} = \mathbf{U}\hat{\mathbf{K}}\mathbf{U}^*$ for some matrix
$\hat{\mathbf{K}} \in \mathbb{C}^{r\times r}$. After some manipulation, it can
be shown that minimizing the empirical error \cref{eqn:EDMD_Frobenius_error}
with respect to $\hat{\mathbf{K}}$ is equivalent to minimizing
\begin{equation}
  \label{eqn:KDMD_Frobenius_error}
  J' = \left\Vert \left(\mathbf{V}^*\mathbf{\Psi}_{\mathbf{Y}}^* - \mathbf{\Sigma}\hat{\mathbf{K}}\mathbf{U}^* \right) \mathbf{A} \right\Vert_F^2.
\end{equation}
Regardless of how the columns of $\mathbf{A}$ are chosen, as long as
$\mathcal{R}(\mathbf{A}) = \mathcal{R}(\mathbf{U})$ the minimum norm solution
for the KDMD matrix is
\begin{equation}
  \label{eqn:KDMD_solution}
  \hat{\mathbf{K}} = \mathbf{\Sigma}^+\mathbf{V}^*\mathbf{\Psi}_{\mathbf{Y}}^*\mathbf{U}
  = \mathbf{\Sigma}^+\mathbf{V}^*\mathbf{\Psi}_{\mathbf{Y}}^*\mathbf{\Psi}_{\mathbf{X}}\mathbf{V}\mathbf{\Sigma}^+.
\end{equation}
Each component in the above KDMD approximation can be found entirely
in terms of inner products in the feature space, enabling the use of a kernel
function to implicitly define the feature space. The two matrices whose entries
are
$\left[\mathbf{K}_{\mathbf{XX}}\right]_{ij} =
\mathbf{\Psi}(\mathbf{x}_i)^*\mathbf{\Psi}(\mathbf{x}_j) = k(\mathbf{x}_i,
\mathbf{x}_j)$ and
$\left[\mathbf{K}_{\mathbf{YX}}\right]_{ij}
=\left[\mathbf{\Psi}_{\mathbf{Y}}^*\mathbf{\Psi}_{\mathbf{X}}\right]_{ij}=
\mathbf{\Psi}(\mathbf{y}_i)^*\mathbf{\Psi}(\mathbf{x}_j) = k(\mathbf{y}_i,
\mathbf{x}_j)$ are computed using the kernel. The Hermitian eigenvalue
decomposition
$\mathbf{K}_{\mathbf{XX}} = \mathbf{V}\mathbf{\Sigma}^2\mathbf{V}^*$ provides
the matrices $\mathbf{V}$ and $\mathbf{\Sigma}$.

It is worth pointing out that the EDMD and KDMD solutions
\cref{eqn:EDMD_solution} and \cref{eqn:KDMD_solution} can be regularized by
truncating the rank $r$ of the SVD
$\mathbf{\Psi}_{\mathbf{X}} = \mathbf{U}\mathbf{\Sigma}\mathbf{V}^*$. In EDMD,
we recognize that
$\mathbf{G} = \frac{1}{M}\mathbf{\Psi}_{\mathbf{X}}\mathbf{\Psi}_{\mathbf{X}}^*
= \frac{1}{M} \mathbf{U}\mathbf{\Sigma}^2\mathbf{U}^*$ is a Hermitian
eigendecomposition. Before finding the pseudoinverse, the rank is truncated to
remove the dyadic components having small singular values.

\subsection{Computing Koopman eigenvalues, eigenfunctions, and modes}
Suppose that $\varphi = \mathbf{\Psi}^*\mathbf{w}$ is an eigenvector of the
Koopman operator in the span of the dictionary with eigenvalue $\mu$. Suppose
also that $\mathbf{w}=\mathbf{U}\hat{\mathbf{w}}$ is in the span of the data in
feature space. From $\mathcal{K}\varphi = \mu\varphi$ it follows that
$\mathbf{\Psi}_{\mathbf{Y}}^*\mathbf{w} = \mu
\mathbf{\Psi}_{\mathbf{X}}^*\mathbf{w}$ by substituting all of the snapshot
pairs. Left-multiplying by $\frac{1}{M}\mathbf{\Psi}_{\mathbf{X}}$ and taking
the pseudoinverse, we obtain
$(\mathbf{G}^+\mathbf{A})\mathbf{w} = \mu(\mathbf{G}^+\mathbf{G})\mathbf{w} =
\mu \mathbf{w}$ where the second equality holds because
$\mathbf{w}\in\mathcal{R}(\mathbf{\Psi}_{\mathbf{X}})$. Therefore, $\mathbf{w}$
is an eigenvector with eigenvalue $\mu$ of the EDMD matrix
\cref{eqn:EDMD_solution}. In terms of the coefficients $\hat{\mathbf{w}}$, we
have
$\mathbf{\Psi}_{\mathbf{Y}}^*\mathbf{U}\hat{\mathbf{w}} = \mu
\mathbf{\Psi}_{\mathbf{X}}^*\mathbf{U}\hat{\mathbf{w}}$, which upon substituting
the definition $\mathbf{\Psi}_{\mathbf{X}}=\mathbf{U}\mathbf{\Sigma}\mathbf{V}$
gives
$\mathbf{\Psi}_{\mathbf{Y}}^*\mathbf{\Psi}_{\mathbf{X}}\mathbf{V}\mathbf{\Sigma}^+\hat{\mathbf{w}}
= \mu \mathbf{V}\mathbf{\Sigma}\hat{\mathbf{w}}$. From the previous statement we
it is evident that
$\mathbf{\Sigma}^+\mathbf{V}^*\mathbf{\Psi}_{\mathbf{Y}}^*\mathbf{\Psi}_{\mathbf{X}}\mathbf{V}\mathbf{\Sigma}^+\hat{\mathbf{w}}
=\hat{\mathbf{K}}\hat{\mathbf{w}}= \mu\hat{\mathbf{w}}$. Hence,
$\hat{\mathbf{w}}$ is an eigenvector of the KDMD matrix \cref{eqn:KDMD_solution}
with eigenvalue $\mu$. Unfortunately, the converses of these statements do not
hold. Nonetheless, approximations of Koopman eigenfunctions,
\begin{equation}
  \label{eqn:Koopman_eigenfunction_EDMD_KDMD}
  \varphi(\mathbf{x}) = \mathbf{\Psi}(\mathbf{x})^*\mathbf{w} = \mathbf{\Psi}(\mathbf{x})^*\mathbf{\Psi}_{\mathbf{X}}\mathbf{V}\mathbf{\Sigma}^+\hat{\mathbf{w}},
\end{equation}
are formed using the right eigenvectors $\mathbf{w}$ and $\hat{\mathbf{w}}$ of
$\mathbf{K}$ and $\hat{\mathbf{K}}$ respectively. In \cref{eqn:Koopman_eigenfunction_EDMD_KDMD} the inner products
$\mathbf{\Psi}(\mathbf{x})^*\mathbf{\Psi}_{\mathbf{X}}$ can be found by
evaluating the kernel function between $\mathbf{x}$ and each point in the
training data $\left\lbrace \mathbf{x}_j \right\rbrace_{j=1}^M$ yielding a
row-vector.

The Koopman modes $\left\lbrace \boldsymbol{\xi}_k \right\rbrace_{k=1}^r$
associated with the full state observable reconstruct the state as a linear
combination of Koopman eigenfunctions. They can be found from the provided
training data using a linear regression process. Let us define the matrices 
\begin{equation}
  \mathbf{\Xi} =
  \begin{bmatrix}
    \boldsymbol{\xi}_1 & \boldsymbol{\xi}_2 & \cdots & \boldsymbol{\xi}_r
  \end{bmatrix}, \qquad
  \mathbf{\Phi}_{\mathbf{X}} =
  \begin{bmatrix}
    \varphi_1(\mathbf{x}_1) & \cdots & \varphi_1(\mathbf{x}_M) \\
    \vdots & \ddots & \vdots \\
    \varphi_r(\mathbf{x}_1) & \cdots & \varphi_r(\mathbf{x}_M)
  \end{bmatrix} =\mathbf{W}_R^T\overline{\mathbf{\Psi}_{\mathbf{X}}} =
  \hat{\mathbf{W}}_R^T\mathbf{\Sigma}\mathbf{V}^T,
\end{equation}
containing the Koopman modes and eigenfunction values at the training points.
In the above, $\mathbf{W}_R$ and $\hat{\mathbf{W}}_R$ are the matrices
whose columns are the right eigenvectors of $\mathbf{K}$ and $\hat{\mathbf{K}}$
respectively. Seeking to linearly reconstruct the state from the eigenfunction
values at each training point, the regression problem,
\begin{equation}
  \label{eqn:Koopman_modes_regression_problem}
  \minimize_{\mathbf{\Xi}\in\mathbb{C}^{n\times r}} \left\Vert \mathbf{X} - \mathbf{\Xi} \mathbf{\Phi}_{\mathbf{X}} \right\Vert_F^2,
\end{equation}
is formulated.
The solution to this standard least squares problem is
\begin{equation}
  \label{eqn:Koopman_modes_least_squares_solution}
  \mathbf{\Xi} = \mathbf{X} \overline{\mathbf{\Psi}_{\mathbf{X}}^+\mathbf{W}_L} =
  \mathbf{X}\overline{\mathbf{V}\mathbf{\Sigma}^+\hat{\mathbf{W}}_L},
\end{equation}
where $\mathbf{W}_L$ and
$\hat{\mathbf{W}}_L$ are the left eigenvector matrices of $\mathbf{K}$ and
$\hat{\mathbf{K}}$ respectively. These matrices must be suitably normalized so
that the left and right eigenvectors form bi-orthonormal sets
$\mathbf{W}_L^*\mathbf{W}_R = \mathbf{I}_r$ and
$\hat{\mathbf{W}}_L^*\hat{\mathbf{W}}_R=\mathbf{I}_r$.
\subsection{Drawbacks of EDMD}
One of the drawbacks of EDMD and KDMD is that the accuracy depends on the chosen
dictionary. For high-dimensional data sets, constructing and evaluating an
explicit dictionary becomes prohibitively expensive. Though the kernel method
allows us to use high-dimensional dictionaries implicitly, the choice of kernel
function significantly impacts the results. In both techniques, higher
resolution is achieved directly by adding more dictionary elements. Therefore,
enormous dictionaries are needed in order to represent complex features. The
shallow representation of features in terms of linear combinations of dictionary
elements means that the effective size of the dictionary must be limited by the
rank of the training data in feature space. As one increases the resolution of
the dictionary, the rank $r$ of the feature space data
$\mathbf{\Psi}_{\mathbf{X}}$ grows and eventually reaches the number of points
$M$ assuming the points are distinct. The number of data points therefore is an
upper bound on the effective number of features we can retain for EDMD or
KDMD. This effective dictionary selection is implicit when we truncate the SVD
of $\mathbf{G}$ or $\mathbf{\Psi}_{\mathbf{X}}$. It is when $r=M$ that we have
retained enough features to memorize the data set up to projection of
$\mathbf{\Psi}_{\mathbf{Y}}$ onto
$\mathcal{R}(\mathbf{\Psi}_{\mathbf{X}})$. Consequently, over-fitting becomes
problematic as we seek dictionaries with high enough resolution to capture
complex features. We illustrate this problem with the following simple example.

\begin{example} \label{ex:EDMD_counterexample}
  Let us consider the linear dynamical system
  \begin{equation}
    \mathbf{x}_{t+1} = \mathbf{f}(\mathbf{x}_{t+1}) =
    \begin{bmatrix}
      1 & 0 \\
      0 & 0.5
    \end{bmatrix} \mathbf{x}_{t+1}
  \end{equation}
  with $\mathbf{x}=\left[x_1, \enskip x_2 \right]^T \in\mathbb{R}^2$. 
  We construct this example to reflect the behavior of EDMD with rich dictionaries
  containing more elements than snapshots. Suppose that we have only two snapshot
  pairs, 
  \begin{equation}
    \mathbf{X} =
    \begin{bmatrix}
      1 & 1 \\
      1 & 0.5
    \end{bmatrix}, \qquad
    \mathbf{Y} =
    \begin{bmatrix}
      1 & 1 \\
      0.5 & 0.25
    \end{bmatrix},
  \end{equation}
  taken by evolving the trajectory two steps from the initial condition $x_0 = \left[ 1, \enskip 1\right]^T$.
  Let us define the following dictionary. Its first two elements are Koopman eigenfunctions whose values are sufficient to describe the full state. 
  In fact, EDMD recovers the original dynamics perfectly from the given snapshots when we take only these first two observables. 
  In this example, we show that by including an extra, unnecessary observable we get a much worse approximation of the dynamics. 
  A third dictionary element which is not an eigenfunction is included in order to demonstrate the effects of an overcomplete dictionary. 
  With these dictionary elements, the data matrices are
  \begin{equation}
    \mathbf{\Psi}(\mathbf{x}) =
    \begin{bmatrix}
      x_1 \\
      x_2 \\
      (x_1)^2 + (x_2)^2
    \end{bmatrix} \quad \Longrightarrow \quad
    \mathbf{\Psi}_{\mathbf{X}} =
    \begin{bmatrix}
      1 & 1 \\
      1 & 0.5 \\
      2 & 1.25
    \end{bmatrix}, \quad
    \mathbf{\Psi}_{\mathbf{Y}} =
    \begin{bmatrix}
      1 & 1 \\
      0.5 & 0.25 \\
      1.25 & 1.0625
    \end{bmatrix}.
  \end{equation}
  Applying \cref{eqn:EDMD_solution} we compute the EDMD matrix and its
  eigendecomposition,
  \begin{equation}
    \mathbf{K} =
    \begin{bmatrix}
      0.9286 & -0.1071 & 0.7321 \\
      -0.2143 & 0.1786 & -0.0536 \\
      0.1429 & 0.2143 & 0.2857
    \end{bmatrix} \quad \Longrightarrow \quad \left\lbrace
      \begin{aligned}
        \mu_1 &= 1.0413\\
        \mu_2 &= 0\\
        \mu_3 & = 0.3515
      \end{aligned} \right. ,
  \end{equation}
  as well as the eigenfunction approximations, 
  \begin{equation}
    \begin{bmatrix}
      \varphi_1(\mathbf{x}) \\
      \varphi_2(\mathbf{x}) \\
      \varphi_3(\mathbf{x})
    \end{bmatrix} = \mathbf{W}_R^T\overline{\mathbf{\Psi}(\mathbf{x})} =
    \begin{bmatrix}
      -0.9627 & 0.2461 & -0.1122 \\
      0.5735 & 0.4915 & 0.5722 \\
      -0.6013 & 0.5722 & 0.5577
    \end{bmatrix}
    \begin{bmatrix}
      x_1 \\
      x_2 \\
      (x_1)^2 + (x_2)^2
    \end{bmatrix}.
  \end{equation}
  It is easy to see that none of the eigenfunctions or eigenvalues are correct for
  the given system even though the learned matrix satisfies
  $\left\Vert \mathbf{\Psi}_{\mathbf{Y}}^* -
    \mathbf{\Psi}_{\mathbf{X}}^*\mathbf{K} \right\Vert_F < 6*10^{-15}$ with 16
  digit precision computed with standard Matlab tools. This shows that even with a
  single additional function in the dictionary, we have severely over-fit the
  data. This is surprising since our original dictionary included two
  eigenfunctions by definition. The nuance comes since EDMD is only guaranteed to
  capture eigenfunctions
  $\varphi(\mathbf{x}) = \mathbf{w}^T\overline{\mathbf{\Psi}(\mathbf{x})}$ where
  $\mathbf{w}$ is in the span of the feature space data
  $\mathcal{R}(\mathbf{\Psi}_{\mathbf{X}})$. In this example, the true
  eigenfunctions do not satisfy this condition; one can check that neither
  $\mathbf{w}=[1,\enskip 0,\enskip 0]^T$ nor
  $\mathbf{w}=[0,\enskip 1,\enskip 0]^T$ is in
  $\mathcal{R}(\mathbf{\Psi}_{\mathbf{X}})$.
\end{example}

\section{Recent approach for dictionary learning}
\label{sec:EDMD_DL}
\begin{figure}[htbp]
  \centering
  \includegraphics[width=0.4\linewidth]{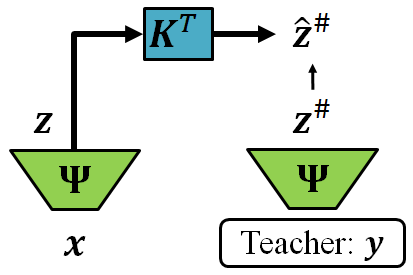}
  \caption{EDMD with dictionary learning architecture}
  \label{fig:EDMD_DL_architecture}
\end{figure}

\Cref{ex:EDMD_counterexample} makes clear the importance of choosing an
appropriate dictionary prior to performing EDMD. In two recent papers
\cite{Yeung2017, Li2017}, the universal function approximation property of
neural networks was used to learn dictionaries for approximating the Koopman
operator. A fixed number of observables making up the dictionary are given by a
neural network $\mathbf{\Psi}(\mathbf{x};\boldsymbol{\theta})\in\mathbb{R}^d$
parameterized by $\boldsymbol{\theta}$. The linear operator
$\mathbf{K}^T\in\mathbb{R}^{d\times d}$ evolving the dictionary function values
one time step into the future is learned simultaneously through minimization of
\begin{equation}
  \label{eqn:EDMD_DL_cost_function}
  J(\mathbf{K}, \boldsymbol{\theta}) = \sum_{i=1}^M \left\Vert \mathbf{\Psi}(\mathbf{y}_i;\boldsymbol{\theta}) - \mathbf{K}^T\mathbf{\Psi}(\mathbf{x}_i;\boldsymbol{\theta}) \right\Vert^2 + \Omega(\mathbf{K}, \boldsymbol{\theta}).
\end{equation}
A schematic of this architecture is depicted
in \cref{fig:EDMD_DL_architecture} where
$\mathbf{z}=\mathbf{\Psi}(\mathbf{x};\boldsymbol{\theta})$ and
$\mathbf{z^{\#}}=\mathbf{\Psi}(\mathbf{y};\boldsymbol{\theta})$ are the
dictionary function values before and after the time increment. The term $\Omega$ is used for regularization and the Tikhonov
regularizer
$\Omega(\mathbf{K}, \boldsymbol{\theta}) = \lambda\left\Vert \mathbf{K}
\right\Vert_F^2$ was used. One notices that as the problem is formulated, the
trivial solution
$\mathbf{\Psi}(\mathbf{x};\boldsymbol{\theta}) \equiv \mathbf{0}_d$ and
$\mathbf{K}=\mathbf{0}_{d\times d}$ is a global minimizer. \cite{Li2017} solves
this problem by fixing some of the dictionary elements to not be
trainable. Since the full state observable is to be linearly reconstructed in
the span of the dictionary elements via the Koopman modes, it is natural to fix
the first $N$ dictionary elements to be $\mathbf{x}$ while learning the
remaining $d-N$ elements through parameterization as a neural network. The
learned dictionary then approximately spans a Koopman invariant subspace
containing the full state observable. Training proceeds by iterating two steps:
(1) Fix $\boldsymbol{\theta}$ and optimize $\mathbf{K}$ by explicit solution of
the least squares problem; Then (2) fix $\mathbf{K}$ and optimize
$\boldsymbol{\theta}$ by gradient descent. The algorithm implemented in
\cite{Li2017} is summarized in \cref{alg:EDMD_DL_training}.

\begin{algorithm}
  \caption{EDMD with dictionary learning \cite{Li2017}}
  \label{alg:EDMD_DL_training}
  \begin{algorithmic}
    \STATE{Initialize $\mathbf{K}$, $\boldsymbol{\theta}$}
    \WHILE{$J(\mathbf{K}, \boldsymbol{\theta}) > \epsilon$}
    \STATE{Tikhonov regularized EDMD: $\mathbf{K} \leftarrow \left(\mathbf{G}(\boldsymbol{\theta}) + \lambda \mathbf{I}_d\right)^{-1}\mathbf{A}(\boldsymbol{\theta})$}
    \STATE{Gradient descent: $\boldsymbol{\theta} \leftarrow \boldsymbol{\theta} - \delta \nabla_{\boldsymbol{\theta}} J(\mathbf{K}, \boldsymbol{\theta})$}
    \ENDWHILE
  \end{algorithmic}
\end{algorithm}

When learning an adaptive dictionary of a fixed size using a neural network (or
other function approximation method), let us consider two objects: the
dictionary space
$\mathcal{S}=\left\lbrace \psi_j(\bullet
  ;\boldsymbol{\theta}):\mathbb{R}^n\rightarrow\mathbb{R} \enskip : \enskip
  \forall \boldsymbol{\theta}\in \Theta, \enskip j=1,\ldots,d \right\rbrace$ is
the set of all functions which can be parameterized by the neural network and
the dictionary
$\mathcal{D}(\boldsymbol{\theta})=\left\lbrace \psi_j(\bullet
  ;\boldsymbol{\theta}):\mathbb{R}^n\rightarrow\mathbb{R} \enskip : \enskip
  j=1,\ldots,d \right\rbrace$ is the $d$ elements of $\mathcal{S}$ fixed by
choosing $\boldsymbol{\theta}$. In EDMD and the dictionary learning approach
just described, the Koopman operator is always approximated on a subspace
$\mathcal{F}_{\mathbf{X}} =\left\lbrace \mathbf{\Psi}^*\mathbf{w}\enskip :
  \enskip \mathbf{w}\in\mathcal{R}(\mathbf{\Psi}_{\mathbf{X}}) \right\rbrace
\subset \mathcal{F}_{\mathcal{D}} = \vspan \mathcal{D}$. As discussed earlier,
the EDMD method always uses $\mathcal{S} = \mathcal{D}$ and the only way to
increase resolution and feature complexity is to grow the dictionary --- leading
to the over-fitting problems illustrated in \cref{ex:EDMD_counterexample}. By
contrast, the dictionary learning approach enables us to keep the size of the
dictionary relatively small while exploring a much larger space
$\mathcal{S}$. In particular, the dictionary size is presumed to be much smaller
than the total number of training data points and probably small enough so that
$\mathcal{F}_{\mathbf{X}} = \mathcal{F}_{\mathcal{D}}$. Otherwise, the number of
functions $d$ learned by the network could be reduced so that this becomes
true. The small dictionary size therefore prevents the method from memorizing the
snapshot pairs without learning true invariant subspaces. This is not done at
the expense of resolution since the allowable complexity of functions in
$\mathcal{S}$ is extremely high.

Deep neural networks are advantageous since they enable highly efficient
representations of certain natural classes of complex features
\cite{Rolnick2017, Lin2016}. In particular, deep neural networks are capable of
learning functions whose values are built by applying many simple operations in
succession. It is shown empirically that this is indeed an important and natural
class of functions since deep neural networks have recently enabled near human
level performance on tasks like image and handwritten digit recognition
\cite{Goodfellow2016}. This proved to be a useful property for dictionary
learning for EDMD since \cite{Yeung2017, Li2017} achieve state of the art
results on examples including the Duffing equation, Kuramoto-Sivashinsky PDE, a
system representing the glycolysis pathway, and power systems.

\section{New approach: deep feature learning using the LRAN}
\begin{figure}[htbp]
  \centering
  \includegraphics[width=0.7\linewidth]{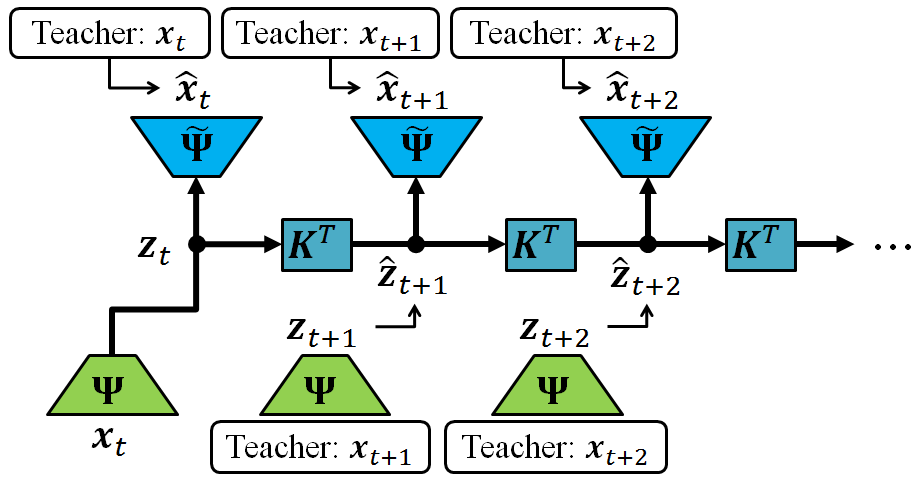}
  \caption{Linearly-Recurrent Autoencoder Network (LRAN) architecture}
  \label{fig:Full_LRAN_architecture}
\end{figure}

By removing the constraint that the full state observable is in the learned
Koopman invariant subspace, one can do even better. This is especially important
for high-dimensional systems where it would be prohibitive to train such a large
neural network-based dictionary with limited training data. Furthermore, it may
simply not be the case that the full state observable lies in a
finite-dimensional Koopman invariant subspace. The method described here is
capable of learning extremely low-dimensional invariant subspaces limited only
by the intrinsic dimensionality of linearly-evolving patterns in the data. A
schematic of our general architecture is presented in
\cref{fig:Full_LRAN_architecture}. The dictionary function values are given by
the output of an encoder neural network
$\mathbf{z}=\mathbf{\Psi}(\mathbf{x};\boldsymbol{\theta}_{enc})$ parameterized
by $\boldsymbol{\theta}_{enc}$. We avoid the trivial solution by nonlinearly
reconstructing an approximation of the full state using a decoder neural network
$\hat{\mathbf{x}} = \tilde{\mathbf{\Psi}}(\mathbf{z};\boldsymbol{\theta}_{dec})$
parameterized by $\boldsymbol{\theta}_{dec}$. The decoder network takes the
place of Koopman modes for reconstructing the full state from eigenfunction
values. However, if Koopman modes are desired it is still possible to compute
them using two methods. The first is to employ the same regression procedure
whose solution is given by \cref{eqn:Koopman_modes_least_squares_solution} to
compute the Koopman modes from the EDMD dictionary provided by the
encoder. Reconstruction using the Koopman modes will certainly achieve lower
accuracy than the nonlinear decoder, but may still provide a useful tool for
feature extraction and visualization. The other option is to employ a linear
decoder network
$\tilde{\mathbf{\Psi}}(\mathbf{z};\boldsymbol{\theta}_{dec}) =
\mathbf{B}(\boldsymbol{\theta}_{dec}) \mathbf{z}$ where
$\mathbf{B}(\boldsymbol{\theta}_{dec})\in\mathbb{R}^{n\times d}$ is a matrix
whose entries are parameterized by $\boldsymbol{\theta}_{dec}$. The advantage of
using a nonlinear decoder network is that the full state observable need not be
in the span of the learned encoder dictionary functions. A nonlinear decoder can
reconstruct more information about the full state from fewer features provided by
the encoder. This enables the dictionary size $d$ to be extremely small --- yet
informative enough to enable nonlinear reconstruction. This is exactly the
principle underlying the success of undercomplete autoencoders for feature
extraction, manifold learning, and dimensionality reduction. Simultaneous
training of the encoder and decoder networks extract rich dictionary elements
which the decoder can use for reconstruction.

The technique includes a linear time evolution process given by the matrix
$\mathbf{K}(\boldsymbol{\theta}_{\mathbf{K}})$ parameterized by
$\boldsymbol{\theta}_{\mathbf{K}}$. This matrix furnishes our approximation of
the Koopman operator on the learned dictionary. Taking the eigendecomposition
$\mathbf{K} = \mathbf{W}_R\mathbf{\Lambda}\mathbf{W}_L^*$ allows us to compute
the Koopman eigenvalues, eigenfunctions, and modes exactly as we would for EDMD
using \cref{eqn:Koopman_eigenfunction_EDMD_KDMD} and
\cref{eqn:Koopman_modes_least_squares_solution}. By training the operator
$\mathbf{K}$ simultaneously with the encoder and decoder networks, the
dictionary of observables learned by the encoder is forced to span a
low-dimensional Koopman invariant subspace which is sufficiently informative to
approximately reconstruct the full state. In many real-world applications, the
scientist has access to data sets consisting of several sequential
snapshots. The LRAN architecture shown in \cref{fig:Full_LRAN_architecture}
takes advantage of longer sequences of snapshots during training. This is
especially important when the system dynamics are highly non-normal. In such
systems, low-amplitude features which could otherwise be ignored for
reconstruction purposes are highly observable and influence the larger amplitude
dynamics several time-steps into the future. One may be able to achieve
reasonable accuracy on snapshot pairs by neglecting some of these low-energy
modes, but accuracy will suffer as more time steps are predicted. Inclusion of
multiple time steps where possible forces the LRAN to incorporate these
dynamically important non-normal features in the dictionary. As we will discuss
later, it is possible to generalize the LRAN architecture to continuous time
systems with irregular sampling intervals and sequence lengths. It is also
possible to restrict the LRAN to the case when only snapshot pairs are
available. Here we consider the case when our data contains equally-spaced
snapshot sequences
$\left\lbrace \mathbf{x}_t,\mathbf{x}_{t+1},\ldots,\mathbf{x}_{t+\mathcal{T}-1}
\right\rbrace$ of length $\mathcal{T}$. The loss function 
\begin{multline}
  \label{eqn:LRAN_general_loss}
 J(\boldsymbol{\theta}_{enc},\boldsymbol{\theta}_{dec},\boldsymbol{\theta}_{\mathbf{K}}) =
  \Expectation_{\mathbf{x}_t,\ldots,\mathbf{x}_{t+\mathcal{T}-1}\sim P_{data}} \frac{1}{1+\beta} \biggl[\sum_{\tau=0}^{\mathcal{T}-1} \frac{\delta^{t}}{N_1(\delta)} \frac{\left\Vert \hat{\mathbf{x}}_{t+\tau} -  \mathbf{x}_{t+\tau} \right\Vert^2}{\left\Vert \mathbf{x}_{t+\tau} \right\Vert^2 + \epsilon_1} \\
  + \beta \sum_{\tau=1}^{\mathcal{T}-1} \frac{\delta^{t-1}}{N_2(\delta)} \frac{\left\Vert \hat{\mathbf{z}}_{t+\tau} -  \mathbf{z}_{t+\tau} \right\Vert^2}{\left\Vert \mathbf{z}_{t+\tau} \right\Vert^2 + \epsilon_2}\biggr] + \Omega(\boldsymbol{\theta}_{enc},\boldsymbol{\theta}_{dec},\boldsymbol{\theta}_{\mathbf{K}})
\end{multline}
is minimized during
training, where $\mathbb{E}$ denotes the expectation over the data distribution.
The encoding, latent state dynamics, and decoding processes are given by
\begin{equation*}
  \mathbf{z}_{t+\tau} = \mathbf{\Psi}(\mathbf{x}_{t+\tau};\boldsymbol{\theta}_{enc}), \quad \hat{\mathbf{z}}_{t+\tau} = \left[\mathbf{K}(\boldsymbol{\theta}_{\mathbf{K}})^{\tau}\right]^T\mathbf{z}_t, \quad \hat{\mathbf{x}}_{t+\tau} = \tilde{\mathbf{\Psi}}(\hat{\mathbf{z}}_{t+\tau};\boldsymbol{\theta}_{dec}),
\end{equation*}
respectively. The regularization term $\Omega$ has been included for generality,
though our numerical experiments show that it was not necessary. Choosing a
small dictionary size $d$ provides sufficient regularization for the
network. The loss function \cref{eqn:LRAN_general_loss} consists of a weighted
average of the reconstruction error and the hidden state time evolution
error. The parameter $\beta$ determines the relative importance of these two
terms. The errors themselves are relative square errors between the predictions
and the ground truth summed over time with a decaying weight $0<\delta\leq
1$. This decaying weight is used to facilitate training by prioritizing short
term predictions. The corresponding normalizing constants,
\begin{equation*}
  N_1(\delta) = \sum_{\tau=0}^{\mathcal{T}-1} \delta^{\tau}, \qquad
  N_2(\delta) = \sum_{\tau=1}^{\mathcal{T}-1} \delta^{\tau-1}
\end{equation*}
ensure that the decay-weighted average is being taken over time. The small constants
$\epsilon_1$ and $\epsilon_2$ are used to avoid division by $0$ in the case that
the ground truth values vanish. The expectation value was estimated empirically using
minibatches consisting of sequences of length $\mathcal{T}$ drawn randomly from
the training data. Stochastic gradient descent with the Adaptive Moment
Estimation (ADAM) method and slowly decaying learning rate was used to
simultaneously optimize the parameters $\boldsymbol{\theta}_{enc}$,
$\boldsymbol{\theta}_{dec}$, and $\boldsymbol{\theta}_{\mathbf{K}}$ in the
open-source software package TensorFlow.

\subsection{Neural network architecture and initialization}
\begin{figure}[htbp]
  \centering
  \includegraphics[width=0.6\linewidth]{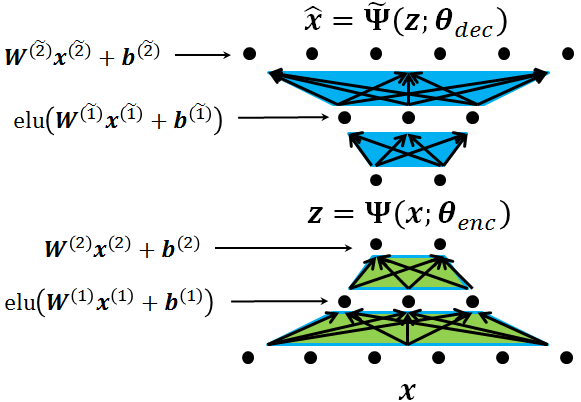}
  \caption{Architecture of the encoder and decoder networks}
  \label{fig:Autoencoder_architecture}
\end{figure}

The encoder and decoder consist of deep neural networks whose schematic is
sketched in \cref{fig:Autoencoder_architecture}. The figure is only a sketch
since many more hidden layers were actually used in the architectures applied to
example problems in this paper. In order to achieve sufficient depth in the
encoder and decoder networks, hidden layers employed exponential linear units or
``elu's'' as the nonlinearity \cite{Clevert2015}. These units mitigate the
problem of vanishing and exploding gradients in deep networks by employing the
identity function for all non-negative arguments. A shifted exponential function
for negative arguments is smoothly matched to the linear section at the origin, giving the activation function
\begin{equation}
  \elu (x) =
  \begin{cases}
    x & \quad x\geq 0 \\
    \exp(x) - 1 & \quad x<0
  \end{cases}.
\end{equation}
This prevents the units from ``dying'' as standard rectified linear
units or ``ReLU's'' do when the arguments are always negative on the data. Furthermore, ``elu's''
have the advantage of being continuously differentiable. This will be a nice
property if we want to approximate a $C^1$ data manifold whose chart map and its
inverse are given by the encoder and decoder. If the maps are differentiable,
then the tangent spaces can be defined as well as push-forward, pull-back, and
connection forms. Hidden layers map the activations $\mathbf{x}^{(l)}$ at layer
$l$ to activations at the next layer $l+1$ given by a linear transformation and subsequent element-wise application of the activation function,
\begin{equation}
  \mathbf{x}^{(l+1)} = \elu \left[ \mathbf{W}^{(l)}(\boldsymbol{\theta}) \mathbf{x}^{(l)} + \mathbf{b}^{(l)}(\boldsymbol{\theta}) \right], \quad
  \mathbf{W}^{(l)}(\boldsymbol{\theta}) \in \mathbb{R}^{n_{l+1}\times n_l}, \quad
  \mathbf{b}^{(l)}(\boldsymbol{\theta}) \in \mathbb{R}^{n_{l+1}}.
\end{equation}
The weight matrices $\mathbf{W}$ and vector biases $\mathbf{b}$ parameterized by
$\boldsymbol{\theta}$ are learned by the network during training. The output
layers $L$ for both the encoder and decoder networks utilize linear transformations without a nonlinear activation function:
\begin{equation}
  \mathbf{x}^{(L)} = \mathbf{W}^{(L-1)}(\boldsymbol{\theta}) \mathbf{x}^{(L-1)} + \mathbf{b}^{(L-1)}(\boldsymbol{\theta}), \quad
  \mathbf{W}^{(L-1)}(\boldsymbol{\theta}) \in \mathbb{R}^{n_{L} \times n_{L-1}}, \quad
  \mathbf{b}^{(L-1)}(\boldsymbol{\theta}) \in \mathbb{R}^{n_{L}},
\end{equation}
where $L=L_{enc}$ or $L=L_{dec}$ is the last layer of the encoder or
decoder with $n_{L_{enc}}=d$ or $n_{L_{dec}}=n$ respectively.
This allows for smooth and consistent
treatment of positive and negative output values without limiting the flow of
gradients back through the network.

The weight matrices were initialized using the Xavier initializer in
Tensorflow. This initialization distributes the entries in $\mathbf{W}^{(l)}$
uniformly over the interval $[-\alpha, \alpha]$ where
$\alpha = \sqrt{6 / (n_l + n_{l+1})}$ in order to keep the scale of gradients
approximately the same in all layers. This initialization method together with
the use of exponential linear units allowed deep networks to be used for the
encoder and decoder. The bias vectors $\mathbf{b}^{(l)}$ were initialized to be
zero. The transition matrix $\mathbf{K}$ was initialized to have diagonal blocks
of the form $\begin{bmatrix} \sigma & \omega \\ -\omega & \sigma \end{bmatrix}$
with eigenvalues $\lambda = \sigma \pm \omega \imath$ equally spaced around the
circle of radius $r = \sqrt{\sigma^2 + \omega^2} = 0.8$. This was done
heuristically to give the initial eigenvalues good coverage of the unit
disc. One could also initialize this matrix using a low-rank DMD matrix.

\subsection{Simple modifications of LRANs}
Several extensions and modifications of the LRAN architecture are possible. Some
simple modifications are discussed here, with several more involved extensions
suggested in the conclusion. In the first extension, we observe that it is easy
learn Koopman eigenfunctions associated with known eigenvalues simply by fixing
the appropriate entries in the matrix $\mathbf{K}$. In particular, if we know
that our system has Koopman eigenvalue $\mu = \sigma + \omega \imath$ then we
may formulate the state transition matrix

\begin{equation}
  \mathbf{K}(\boldsymbol{\theta}) = \begin{bmatrix}
    \begin{bmatrix}
      \sigma & \omega \\
      -\omega & \sigma
    \end{bmatrix} & \mathbf{0}_{2\times (d-2)} \\
    \mathbf{0}_{(d-2) \times 2} & \tilde{\mathbf{K}}(\boldsymbol{\theta})
  \end{bmatrix}.
\end{equation}

\noindent In the above, the known eigenvalue is fixed and only the entries of
$\tilde{\mathbf{K}}$ are allowed to be trained. If more eigenvalues are known,
we simply fix the additional entries of $\mathbf{K}$ in the same way. The case
where some eigenvalues are known is particularly interesting because in certain
cases, eigenvalues of the linearized system are Koopman eigenvalues whose
eigenfunctions have useful properties. Supposing the autonomous system under
consideration has a fixed point with all eigenvalues $\mu_i$, $i=1,\ldots,n$
inside the unit circle, the Hartman-Grobman theorem establishes a topological
conjugacy to a linear system with the same eigenvalues in a small neighborhood
$\mathcal{U}$ of the fixed point. One easily checks that the coordinate
transformations $h_i : \mathcal{M} \cap \mathcal{U} \rightarrow \mathbb{R}$,
$i=1,\ldots,n$ establishing this conjugacy are Koopman eigenfunctions restricted
to the neighborhood. Composing them with the flow map allows us to extend the
eigenfunctions to the entire basin of attraction by defining
$\varphi_i(\mathbf{x}) = \mu_i^{-\tau(\mathbf{x})}
h_i\left(\mathbf{f}^{\tau(\mathbf{x})} (\mathbf{x}) \right)$ where
$\tau(\mathbf{x})$ is the smallest integer $\tau$ such that
$\mathbf{f}^{\tau}(\mathbf{x})\in\mathcal{U}$. These eigenfunctions extend the
topological conjugacy by parameterizing the basin. Similar results hold for
stable limit cycles and tori \cite{Mezic2017}. This is nice because we can often
find eigenvalues at fixed points explicitly by linearization. Choosing to fix
these eigenvalues in the $\mathbf{K}$ matrix forces the LRAN to learn
corresponding eigenfunctions parameterizing the basin of attraction. It is also
easy to include a set of observables explicitly by appending them to the encoder
function
$\mathbf{\Psi}(\mathbf{x};\boldsymbol{\theta}_{enc}) = \left[
  \mathbf{\Psi}_{fixed}(\mathbf{x})^T,\enskip
  \tilde{\mathbf{\Psi}}(\mathbf{x};\boldsymbol{\theta}_{enc})^T \right]^T$ so
that only the functions $\tilde{\mathbf{\Psi}}$ are learned by the network. This
may be useful if we want to accurately reconstruct some observables
$\mathbf{\Psi}_{fixed}$ linearly using Koopman modes.

The LRAN architecture and loss function \cref{eqn:LRAN_general_loss} may be
further generalized to non-uniform sampling of continuous-time systems. In this
case, we consider $T$ sequential snapshots \\
$\left\lbrace \mathbf{x}(t_0), \mathbf{x}(t_1), \ldots, \mathbf{x}(t_{T-1})
\right\rbrace$ where the times $t_0, t_1,\ldots,t_{T-1}$ are not necessarily
evenly spaced. In the continuous time case, we have a Koopman operator semigroup
$\mathcal{K}_{t+s}=\mathcal{K}_{t}\mathcal{K}_{s}$ defined as
$\mathcal{K}_t\psi(\mathbf{x}) = \psi(\mathbf{\varPhi}_t(\mathbf{x}))$ and
generated by the operator
$\mathcal{K}\psi(\mathbf{x}) = \dot{\psi}(\mathbf{x}) =
\left(\nabla_{\mathbf{x}} \psi(\mathbf{x})\right)\mathbf{f}(\mathbf{x})$ where
the dynamics are given by $\dot{\mathbf{x}}=\mathbf{f}(\mathbf{x})$ and
$\mathbf{\varPhi}_t$ is the time $t$ flow map. The generator $\mathcal{K}$ is
clearly a linear operator which we can approximate on our dictionary of
observables with a matrix $\mathbf{K}$. By integrating, we can approximate
elements from the semigroup $\mathcal{K}_t$ using the matrices
$\mathbf{K}_t = \exp{(\mathbf{K}t)}$ on the dictionary. Finally, in order to
formulate the analogous loss function, we might utilize continuously decaying weights
\begin{equation}
  \rho_1 (t) = \frac{\delta^t}{\sum_{k=0}^{\mathcal{T}-1}\delta^{t_k}}, \qquad
  \rho_2 (t) = \frac{\delta^t}{\sum_{k=1}^{\mathcal{T}-1}\delta^{t_k}},
\end{equation}
normalized so that they sum to $1$ for the given sampling times. 
Neural networks can be used for the encoder and decoder
together with the loss function
\begin{multline}
\label{eqn:LRAN_continuous_loss}
J(\boldsymbol{\theta}_{enc},\boldsymbol{\theta}_{dec},\boldsymbol{\theta}_{\mathbf{K}}) =
\Expectation_{\mathbf{x}(t_0),\ldots,\mathbf{x}(t_{\mathcal{T}-1})\sim P_{data}} \frac{1}{1+\beta} \biggl[\sum_{k=0}^{\mathcal{T}-1} \rho_1(t_k) \frac{\left\Vert \hat{\mathbf{x}}(t_k) -  \mathbf{x}(t_k) \right\Vert^2}{\left\Vert \mathbf{x}(t_k) \right\Vert^2 + \epsilon_1} \\
+ \beta \sum_{k=1}^{\mathcal{T}-1} \rho_2(t_k) \frac{\left\Vert \hat{\mathbf{z}}(t_k) -  \mathbf{z}(t_k) \right\Vert^2}{\left\Vert \mathbf{z}(t_k) \right\Vert^2 + \epsilon_2}\biggr] + \Omega(\boldsymbol{\theta}_{enc},\boldsymbol{\theta}_{dec},\boldsymbol{\theta}_{\mathbf{K}})
\end{multline}
to be minimized during training. In this case, the dynamics evolve the observables linearly in continuous time, so we let
\begin{equation*}
  \mathbf{z}(t_k) = \mathbf{\Psi}(\mathbf{x}(t_k);\boldsymbol{\theta}_{enc}), \quad \hat{\mathbf{z}}(t_k) = \exp\left[\mathbf{K}(\boldsymbol{\theta}_{\mathbf{K}}) (t_k-t_0)\right]^T\mathbf{z}(t_0), \quad \hat{\mathbf{x}}(t_k) = \tilde{\mathbf{\Psi}}(\hat{\mathbf{z}}(t_k);\boldsymbol{\theta}_{dec}).
\end{equation*}

\noindent This loss function can be evaluated on the training data and minimized
in essentially the same way as \cref{eqn:LRAN_general_loss}. The only difference
is that we are discovering a matrix approximation to the generator of the
Koopman semigroup. We will not explore irregularly sampled continuous time
systems further in this paper, leaving it as a subject for future work.

\begin{figure}[htbp]
  \centering
  \includegraphics[width=0.4\linewidth]{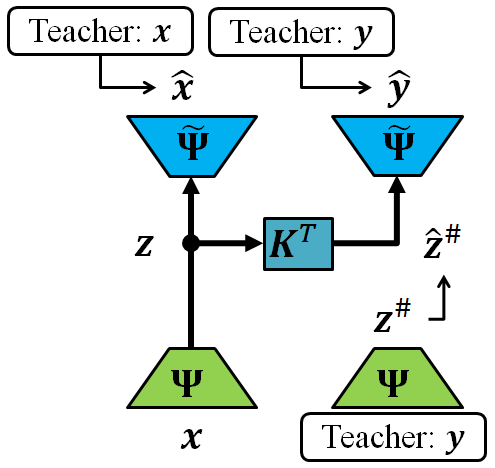}
  \caption{LRAN architecture restricted to snapshot pairs}
  \label{fig:Sing_Step_LRAN_architecture}
\end{figure}

We briefly remark that the general LRAN architecture can be restricted to the
case of snapshot pairs as shown in \cref{fig:Sing_Step_LRAN_architecture}. In
this special case, training might be accelerated using a technique similar to
\cref{alg:EDMD_DL_training}. During the initial stage of training, it may be
beneficial to periodically re-initialize the $\mathbf{K}$ matrix with its EDMD
approximation using the current dictionary functions and a subset of the
training data. This might provide a more suitable initial condition for the
matrix as well as accelerate the training process. However, this update for
$\mathbf{K}$ is not consistent with all the terms in the loss function $J$ since
it does not account for reconstruction errors. Therefore, the final stages of
training must always proceed by gradient descent on the complete loss function.

Finally, we remark that the LRAN architecture sacrifices linear reconstruction
using Koopman modes for nonlinear reconstruction using a decoder neural network
in order to learn ultra-low dimensional Koopman invariant
subspaces. Interestingly, this formulation allows the LRAN to parameterize the
Nonlinear Normal Modes (NNMs) frequently encountered in structural
dynamics. These modes are two-dimensional, periodic invariant manifolds
containing a fixed point of a Hamiltonian system lacking internal
resonances. Therefore, if $\mu = \omega \imath$ and
$\overline{\mu}=-\omega \imath$ are a complex conjugate pair of pure imaginary
eigenvalues of $\mathbf{K}$ with corresponding left eigenvectors $\mathbf{w}_L$
and $\overline{\mathbf{w}_L}$ then a NNM is parameterized as follows:
\begin{equation}
  \mathbf{x}(\alpha) = \tilde{\mathbf{\Psi}}(\alpha\mathbf{w}_L+\overline{\alpha\mathbf{w}_L};\boldsymbol{\theta}_{dec}) = \tilde{\mathbf{\Psi}}\left(2\Re(\alpha)\Re(\mathbf{w}_L)-2\Im(\alpha)\Im(\mathbf{w}_L);\boldsymbol{\theta}_{dec}\right).
\end{equation}
The global coordinates on the manifold are
$\left(\Re(\alpha), \Im(\alpha_2)\right)$. Coordinate projection of the full
state onto the NNM,
\begin{equation}
  \begin{bmatrix}
    \Re(\alpha) \\
    \Im(\alpha)
  \end{bmatrix} =
  \begin{bmatrix}
    \Re(\mathbf{w}_R)^T \\
    \Im(\mathbf{w}_R)^T
  \end{bmatrix} \mathbf{\Psi}(\mathbf{x};\boldsymbol{\theta}_{enc}),
\end{equation}
is accomplished by employing the encoder network and the right eigenvector $\mathbf{w}_R$ corresponding to eigenvalue $\mu$.
These coordinates are the real and imaginary parts of the associated
Koopman eigenfunction $\alpha = \varphi(\mathbf{x})$. The NNM has angular
frequency $\angle(\omega\imath)/\Delta t$ where $\Delta t$ is the sampling
interval between the snapshots in the case of the discrete time LRAN.

We may further generalize the notion of NNMs by considering the Koopman mode
expansion of the real-valued observable vector $\mathbf{\Psi}$ making up our
dictionary. In this particular case, the associated Koopman modes are the
complex conjugate left eigenvectors of $\mathbf{K}$. They allow exact
reconstruction and prediction using the decomposition
\begin{equation}
  \mathbf{z}_{t} = \sum_{j=1}^{r}\overline{\mathbf{w}_{L,j}} \mu_j^t \left(\mathbf{w}_{R,j}^T\mathbf{z}_0\right) = \sum_{j=1}^{r}\overline{\mathbf{w}_{L,j}} \mu_j^t \varphi_j(\mathbf{x}_0),
\end{equation}
assuming a Koopman invariant subspace has been learned that contains the full state observable. 
Reconstructing and predicting with the decoder instead, we have
\begin{equation}
  \mathbf{x}_{t} = \tilde{\mathbf{\Psi}}\left[ \sum_{j=1}^{r}\overline{\mathbf{w}_{L,j}} \mu_j^t \varphi_j(\mathbf{x}_0);\enskip \boldsymbol{\theta}_{dec}\right].
\end{equation}
Therefore, each invariant subspace of $\mathbf{K}$ given by its left
eigenvectors corresponds to an invariant manifold in the $n$-dimensional phase
space. These manifolds have global charts whose coordinate projections are given
by the Koopman eigenfunctions
$\varphi_j(\mathbf{x}) =
\mathbf{w}_{R,j}^T\mathbf{\Psi}(\mathbf{x};\boldsymbol{\theta}_{enc})$. The
dynamics on these manifolds is incredibly simple and entails repeated
multiplication of the coordinates by the eigenvalues. Generalized eigenspaces
may also be considered in the natural way by using the Jordan normal form of
$\mathbf{K}$ instead of its eigendecomposition in the above arguments. The only
necessary change is in the evolution equations, where instead of taking powers
of $\mathbf{\Lambda} = \diag\left\lbrace\mu_1,\ldots,\mu_r\right\rbrace$, we
take powers of $\mathbf{J}$, the matrix consisting of Jordan blocks
\cite{Mezic2017}. Future work might use a variational autoencoder (VAE) formulation \cite{Goodfellow2016,Mescheder2017,Makhzani2015,Krishnan2017} where a given distribution is imposed on the latent state in order to facilitate sampling.

\section{EDMD-based model reduction as a shallow autoencoder}
In this section we examine how the EDMD method might be used to construct
low-dimensional Koopman invariant subspaces while still allowing for accurate
reconstructions and predictions of the full state. The idea is to find a reduced
order model of the large linear system identified by EDMD in the space of
observables. This method is sometimes called overspecification \cite{Rowley2017}
and essentially determines an encoder function into an appropriate small set of
features. From this reduced set of features, we then employ nonlinear
reconstruction of the full state through a learned decoder
function. Introduction of the nonlinear decoder should allow for
lower-dimensional models to be identified which are still able to make accurate
predictions. The proposed framework therefore constructs a kind of autoencoder
where encoded features evolve with linear time invariant dynamics. The encoder
functions are found explicitly as linear combinations of EDMD observables and
are therefore analogous to a shallow neural network with a single hidden
layer. The nonlinear decoder function is also found explicitly through a
regression process involving linear combinations of basis functions.

We remark that this approach differs from training an LRAN by minimization of
\cref{eqn:LRAN_general_loss} in two important ways. First, the EDMD-based model
reduction and reconstruction processes are performed sequentially; thus, the parts 
are not simultaneously optimized as in the LRAN. The LRAN is advantageous since
we only learn to encode observables which the decoder can successfully use for
reconstruction. There are no such guarantees here. Second, the EDMD dictionary
remains fixed albeit overspecified whereas the LRAN explicitly learns an
appropriate dictionary. Therefore, the EDMD shallow autoencoder framework will
still suffer from the overfitting problem illustrated in
\cref{ex:EDMD_counterexample}. If the EDMD-identified matrix $\mathbf{K}$ does
not correctly represent the dynamics on a Koopman invariant subspace, then any
reduced order models derived from it cannot be expected to be accurately relect
the dynamics either. Nonetheless, in many cases, this method could provide a
less computationally expensive alternative to training a LRAN which retains some
of the benefits owing to nonlinear reconstruction.

Dimensionality reduction is achieved by first performing EDMD with a large
dictionary, then projecting the linear dynamics onto a low-dimensional
subspace. A naive approach would be to simply project the large feature space
system onto the most energetic POD modes --- equivalent to low-rank truncation
of the SVD
$\mathbf{\Psi}_{\mathbf{X}}=\mathbf{U}\mathbf{\Sigma}\mathbf{V}^*$. While
effective for normal systems with a few dominant modes, this approach yields
very poor predictions in non-normal systems since low amplitude modes with large
impact on the dynamics would be excluded from the model. One method which
resolves this issue is balanced truncation of the identified feature space
system. Such an idea is suggested in \cite{Rowley2017} for reducing the system
identified by linear DMD. Drawing from the model reduction procedure for
snapshot-based realizations developed in \cite{Luchtenburg2011}, we will
construct a balanced reduced order model for the system identified using EDMD or
KDMD. In the formulation of EDMD that led to the kernel method, an
approximation of the Koopman operator,
\begin{equation}
  \mathcal{K}\mathbf{\Psi}(\mathbf{x})^*\mathbf{a} = \mathbf{\Psi}(\mathbf{x})^*\mathbf{U}\hat{\mathbf{K}}\mathbf{U}^*\mathbf{a} +\mathbf{r}(\mathbf{x}), \quad \forall \mathbf{a}\in\mathcal{R}(\mathbf{U}),
\end{equation}
was obtained. The approximation allows us to model the dynamics of a vector of observables,
\begin{equation}
  \mathbf{\Psi}_{\mathbf{U}}(\mathbf{x})=\mathbf{U}^*\mathbf{\Psi}(\mathbf{x})=\mathbf{\Sigma}^+\mathbf{V}^*\mathbf{\Psi}_{\mathbf{X}}^*\mathbf{\Psi}(\mathbf{x}),
\end{equation}
with the linear input-output system
\begin{equation}
  \label{eqn:feature_state_space_system}
  \begin{aligned}
    \mathbf{\Psi}_{\mathbf{U}}(\mathbf{x}_{t+1}) &=  \hat{\mathbf{K}}^*\mathbf{\Psi}_{\mathbf{U}}(\mathbf{x}_{t}) + \dfrac{1}{\sqrt{M}}\mathbf{\Sigma}\mathbf{u}_t \\
    \mathbf{x}_t &= \mathbf{C}\mathbf{\Psi}_{\mathbf{U}}(\mathbf{x}_{t})
  \end{aligned},
\end{equation}
where $\hat{\mathbf{K}}$ is the matrix \cref{eqn:KDMD_solution}
identified by EDMD or KDMD. The input $\mathbf{u}_t$ is provided in order to
equate varying initial conditions $\mathbf{\Psi}_{\mathbf{U}}(\mathbf{x}_{0})$
with impulse responses of \cref{eqn:feature_state_space_system}. Since the input
is used to set the initial condition, we choose to scale each component by its
singular value to reflect the covariance
\begin{equation}
  \Expectation_{\mathbf{x}\sim P_{data}}\left[\mathbf{\Psi}_{\mathbf{U}}(\mathbf{x})\mathbf{\Psi}_{\mathbf{U}}(\mathbf{x})^*\right] \approx \frac{1}{M} \mathbf{U}^*\mathbf{\Psi}_{\mathbf{X}}\mathbf{\Psi}_{\mathbf{X}}^*\mathbf{U} = \frac{1}{M}\mathbf{\Sigma}^2=\Expectation_{\mathbf{u}\sim \mathcal{N}(\mathbf{0},\mathbf{I}_r)}\left[\dfrac{1}{M}\mathbf{\Sigma}\mathbf{u}\mathbf{u}^*\mathbf{\Sigma}^*\right],
\end{equation}
in the observed data. Therefore, initializing the system using impulse responses
$\mathbf{u}_0 = \mathbf{e}_j$, $j=1,\ldots,r$ from the $\sigma$-points of the
distribution $\mathbf{u}\sim \mathcal{N}(\mathbf{0},\mathbf{I}_r)$ ensures that
the correct empirical covariances are obtained. The output matrix,
\begin{equation}
  \mathbf{C} = \mathbf{X}\mathbf{V}\mathbf{\Sigma}^+,
\end{equation}
is used to linearly reconstruct the full state observable from the complete set
of features. It is found using linear regression similar to the Koopman modes
\cref{eqn:Koopman_modes_least_squares_solution}. The low-dimensional set of observables making up the encoder will be
found using a balanced reduced order model of
\cref{eqn:feature_state_space_system}.

\subsection{Balanced Model Reduction}
Balanced truncation \cite{Moore1981} is a projection-based model reduction
techniqe that attempts to retain a subspace in which
\cref{eqn:feature_state_space_system} is both maximally observable and
controllable. While these notions generally do not coincide in the original
space, remarkably it is possible to find a left-invertible linear transformation
$\mathbf{\Psi}_{\mathbf{U}}(\mathbf{x}) = \mathbf{T}\mathbf{z}$ under which
these properties are balanced. This so called balancing transformation of the
learning subspace simultaneously diagonalizes the observability and
controllability Gramians. Therefore, the most observable states are also the
most controllable and vice versa. The reduced order model is formed by
truncating the least observable/controllable states of the transformed
system. If the discrete time observability Gramian $\mathbf{W}_o$ and
controllability Gramian $\mathbf{W}_c$ are given by
\begin{equation}
  \mathbf{W}_o = \sum_{t=0}^{\infty} (\hat{\mathbf{K}})^t\mathbf{C}^*\mathbf{C}(\hat{\mathbf{K}}^*)^t, \qquad
  \mathbf{W}_c = \frac{1}{M} \sum_{t=0}^{\infty} (\hat{\mathbf{K}}^*)^t\mathbf{\Sigma}^2(\hat{\mathbf{K}})^t,
\end{equation}
then the Gramians transform according to 
\begin{equation}
  \mathbf{W}_o \mapsto \mathbf{T}^*\mathbf{W}_o\mathbf{T}, \qquad
  \mathbf{W}_c \mapsto \mathbf{T}_L^+\mathbf{W}_c(\mathbf{T}_L^+)^*
\end{equation}
under the change of coordinates.
In the above, $\mathbf{S}^*=\mathbf{T}_L^+$ is the left pseudoinverse
satisfying $\mathbf{S}^*\mathbf{T}=\mathbf{I}_d$ where $d$ is the rank of
$\mathbf{T}$ and $\mathbf{T}\mathbf{S}^*=\mathbf{P}_{\mathbf{T}}$ is a (not
necessarily orthogonal) projection operator onto $\mathcal{R}(\mathbf{T})$.

Since the Gramians are Hermitian positive semidefinite, they can be written as
$\mathbf{W}_o = \mathbf{A}\mathbf{A}^*$, $\mathbf{W}_c = \mathbf{B}\mathbf{B}^*$
for some not necessarily unique matrices
$\mathbf{A},\mathbf{B}\in\mathbb{R}^{r\times r}$. Forming an economy sized
singular value decomposition
$\mathbf{H}=\mathbf{A}^*\mathbf{B}=\mathbf{U}_{H}\mathbf{\Sigma}_H\mathbf{V}_H$
allows us to construct the transformations
\begin{equation}
  \mathbf{T} = \mathbf{B}\mathbf{V}_H(\mathbf{\Sigma}_H^+)^{1/2}, \qquad
  \mathbf{S} = \mathbf{A}\mathbf{U}_H(\mathbf{\Sigma}_H^+)^{1/2}.
\end{equation}
Using this construction, it is easy to check that the resulting
transformation simultaneously diagonalizes the Gramians:
\begin{equation}
  \mathbf{T}^*\mathbf{W}_o\mathbf{T} = \mathbf{S}^*\mathbf{W}_c \mathbf{S} = \mathbf{\Sigma}_H.
\end{equation}
 \noindent Entries of the diagonal matrix $\mathbf{\Sigma}_H$ are called the
 Hankel singular values. The columns of $\mathbf{T}$ and $\mathbf{S}$ are called
 the ``balancing modes'' and ``adjoint modes'' respectively. The balancing modes
 span the subspace where \cref{eqn:feature_state_space_system} is both
 observable and controllable while the adjoint modes furnish the projected
 coefficients of the state
 $\mathbf{S}^* \mathbf{\Psi}_{\mathbf{U}}(\mathbf{x}) = \mathbf{z}$ onto the
 space where these properties are balanced. The corresponding Hankel singular
 values quantify the observability/controllability of the states making up
 $\mathbf{z}$. Therefore, a reduced order model which is provably close to
 optimal truncation in the $\mathcal{H}_{\infty}$ norm is formed by rank-$d$ truncation of
 the SVD, retaining only the first $d$ balancing and adjoint modes
 $\mathbf{T}_d$ and $\mathbf{S}_d$ \cite{Dullerud2000}. The reduced state space
 system modeling the dynamics of \cref{eqn:feature_state_space_system} is given by
 \begin{equation}
   \label{eqn:reduced_feature_state_space_system}
   \begin{aligned}
     \mathbf{z}_{t+1} &= \mathbf{S}_d^*\hat{\mathbf{K}}^*\mathbf{T}_d\mathbf{z}_t + \dfrac{1}{\sqrt{M}}\mathbf{S}_d^*\mathbf{\Sigma}\mathbf{u}_t \\
     \mathbf{x}_{t} &\approx \mathbf{C}\mathbf{T}_d\mathbf{z}_t
   \end{aligned}, \quad \mbox{where} \quad \mathbf{z} \triangleq \mathbf{S}_d^*\mathbf{\Psi}_{\mathbf{U}}(\mathbf{x}).
 \end{equation}
Therefore, the reduced dictionary of observables is given by the
components of
$\mathbf{S}_d^*\mathbf{\Psi}_{\mathbf{U}}:\mathbb{R}^n \rightarrow \mathbb{R}^d$
and the corresponding EDMD approximation of the Koopman operator on its span is
given by $\mathbf{T}_d^*\hat{\mathbf{K}}\mathbf{S}_d$. This set of features is
highly observable in that their dynamics strongly influences the full state
reconstruction over time through $\mathbf{C}$. And highly controllable in that
the features are excited by typical state configurations through
$\mathbf{\Sigma}$. The notion of feature excitation corresponding to
controllability will be made clear in the next section.

\subsection{Finite-horizon Gramians and Balanced POD}
Typically one would find the infinite horizon Gramians for an overspecified
Hurwitz EDMD system \cref{eqn:feature_state_space_system} by solving the
Lyapunov equations
\begin{equation}
  \label{eqn:Lyapunov_equations}
  \hat{\mathbf{K}}\mathbf{W}_o \hat{\mathbf{K}}^* - \mathbf{W}_o
  + \mathbf{C}^*\mathbf{C} = 0,
  \qquad \hat{\mathbf{K}}^*\mathbf{W}_c \hat{\mathbf{K}} - \mathbf{W}_c
  + \mathbf{\Sigma}^2 = 0.
\end{equation}

In the case of neutrally stable or unstable systems, unique positive definite
solutions do not exist and one must use generalized Gramians
\cite{Zhou1999}. When used to form balanced reduced order models, this will
always result in the unstable and neutrally stable modes being chosen before the
stable modes. This could be problematic for our intended framework since EDMD
can identify many spurious and sometimes unstable eigenvalues corresponding to
noisy low-amplitude fluctuations. While these noisy modes remain insignificant
over finite times of interest, they will dominate EDMD-based predictions over
long times. Therefore it makes sense to consider the dominant modes identified
by EDMD over a finite time interval of interest. Using finite horizon Gramians
reduces the effect of spurious modes on the reduced order model, making
it more consistent with the data. The time horizon can be chosen to reflect a
desired future prediction time or the number of sequential snapshots in the
training data.

The method of Balanced Proper Orthogonal Decomposition or BPOD \cite{Rowley2005}
allows us to find balancing and adjoint modes of the finite horizon system. In
BPOD, we observe that the finite-horizon Gramians are empirical covariance
matrices formed by evolving the dynamics from impulsive initial conditions for
time $\mathcal{T}$. This gives the specific form for matrices 
\begin{equation}
  \mathbf{A} =
  \begin{bmatrix}
    \mathbf{C}^* & \hat{\mathbf{K}}\mathbf{C}^* & \cdots & (\hat{\mathbf{K}})^{\mathcal{T}}\mathbf{C}^*
  \end{bmatrix}, \qquad
  \mathbf{B} =
  \dfrac{1}{\sqrt{M}}\begin{bmatrix}
    \mathbf{\Sigma} & \hat{\mathbf{K}}^*\mathbf{\Sigma} & \cdots & (\hat{\mathbf{K}}^*)^{\mathcal{T}}\mathbf{\Sigma}
  \end{bmatrix},
\end{equation}
allowing for computation of the balancing and adjoint modes
without ever forming the Gramians. This is known as the method of snapshots.
Since the output dimension is large, we consider its projection onto
the most energetic modes. These are identified by forming the economy sized SVD
of the impulse responses
$\mathbf{C}\mathbf{B}=\mathbf{U}_{OP}\mathbf{\Sigma}_{OP}\mathbf{V}_{OP}^*$. Projecting
the output allows us to form the elements of
\begin{equation}
  \mathbf{A}_{OP} =
  \begin{bmatrix}
    \mathbf{C}^*\mathbf{U}_{OP} & \hat{\mathbf{K}}\mathbf{C}^*\mathbf{U}_{OP} & \cdots & (\hat{\mathbf{K}})^{\mathcal{T}}\mathbf{C}^*\mathbf{U}_{OP}
  \end{bmatrix},
\end{equation}
from fewer initial
conditions than $\mathbf{A}$. In particular, the initial conditions are the first few columns of $\mathbf{U}_{OP}$ with the
largest singular values \cite{Rowley2005}.

Observe that the unit impulses place the initial conditions precisely at the
$\sigma$-points of the data-distribution in features space. If this distribution
is Gaussian, then the empirical expectations obtained by evolving the linear
system agree with the true expectations taken over the entire data
distribution. Therefore, the finite horizon controllability Gramian corresponds
to the covariance matrix taken over all time $\mathcal{T}$ trajectories starting
at initial data points coming from a Gaussian distribution in feature
space. Consequently, controllability in this case corresponds exactly with
feature variance or expected square amplitude over time.

We remark that in the infinite-horizon limit $\mathcal{T}\rightarrow\infty$,
BPOD converges on a transformation which balances the generalized Gramians
introduced in \cite{Zhou1999}. Application of BPOD to unstable systems is
discussed in \cite{Flinois2015} which provides justification for the approach.

Another option to avoid spurious modes from corrupting the long-time dynamics is
to consider pre-selection of EDMD modes which are nearly Koopman invariant. The
development of such an accuracy criterion for selecting modes is the subject of
a forthcoming paper by H. Zhang and C. W. Rowley. One may then apply balanced
model reduction to the feature space system containing only the most accurate
modes.

\subsection{Nonlinear reconstruction}
In truncating the system, we determined a small subspace of observables whose
values evolve linearly in time and are both highly observable and
controllable. However, most of the less observable low-amplitude modes are removed. 
The linear reconstruction only allows us to represent data on a
low-dimensional subspace of $\mathbb{R}^n$. While projection onto this subspace
aims to explain most of the data variance and dynamics, it may be the case that
the data lies near a curved manifold not fully contained in the subspace. The
neglected modes contribute to this additional complexity in the shape of the
data in $\mathbb{R}^n$. Nonlinear reconstruction of the full state can help
account for the complex shape of the data and for neglected modes enslaved to
the ones retained.

We consider the regression problem involved in reconstructing the full state
$\mathbf{x}$ from a small set of EDMD observables $\mathbf{z}$. Because the
previously obtained solution to the EDMD balanced model reduction problem
\cref{eqn:reduced_feature_state_space_system} employs linear reconstruction
through matrix $\mathbf{C}\mathbf{T}_d$, we expect nonlinearities in the
reconstruction to be small with most of the variance being accounted for by
linear terms. Therefore, the regression model,
\begin{equation}
  \label{eqn:reconstruction_regression_model}
  \mathbf{x} = \mathbf{C}_1 \mathbf{z} + \mathbf{C}_2 \mathbf{\Psi}(\mathbf{z}) + \mathbf{e},
\end{equation}
is formulated based on
\cite{Espinoza2005,Xu2009} to include linear and nonlinear components.
In the above, $\mathbf{\Psi}:\mathbb{C}^d \rightarrow \mathcal{H}$ is
a nonlinear feature map into reproducing kernel Hilbert space $\mathcal{H}$ and
$\mathbf{C}_1:\mathbb{C}^d\rightarrow\mathbb{C}^n$ and
$\mathbf{C}_2:\mathcal{H}\rightarrow\mathbb{C}^n$ are linear operators. These operators are found by solving the $l^2$ regularized optimization problem,
\begin{equation}
  \label{eqn:multikernel_regression_loss_1}
  \minimize_{\mathbf{C}_1, \mathbf{C}_2} J = \left\Vert \mathbf{X}^*
    - \mathbf{Z}^*\mathbf{C}_1^* - \mathbf{\Psi}_{\mathbf{Z}}^*\mathbf{C}_2^* \right\Vert_F^2
  + \gamma \Tr (\mathbf{C}_2\mathbf{C}_2^*),
  \qquad \gamma \geq 0,
\end{equation}
involving the empirical square error on the training data
$\left\lbrace(\mathbf{z}_j, \mathbf{x}_j)\right\rbrace_{j=1}^M$ arranged into
columns of the matrices
$\mathbf{Z}=\begin{bmatrix} \mathbf{z}_1 & \cdots & \mathbf{z}_M \end{bmatrix}$
and
$\mathbf{X}=\begin{bmatrix} \mathbf{x}_1 & \cdots &
  \mathbf{x}_M \end{bmatrix}$.
The regularization penalty is placed only on the coefficients of nonlinear terms to control over-fitting while the linear term, which we expect to dominate, is not penalized.
The operator
$\mathbf{\Psi}_{\mathbf{Z}}:\mathbb{C}^M\rightarrow\mathcal{H}$ forms linear
combinations of the data in feature space
$\mathbf{v} \mapsto v_1\mathbf{\Psi}(\mathbf{z}_1) + \cdots +
v_M\mathbf{\Psi}(\mathbf{z}_M)$. Since $\mathbf{Z}$ and
$\mathbf{\Psi}_{\mathbf{Z}}$ are operators with finite ranks $r_1$ and
$r_2\leq M$, we may consider their economy sized singular value decompositions:
$\mathbf{Z} = \mathbf{U}_1\mathbf{\Sigma}_1\mathbf{V}_1^*$ and
$\mathbf{\Psi}_{\mathbf{Z}} =
\mathbf{U}_2\mathbf{\Sigma}_2\mathbf{V}_2^*$. Observe that it is impossible to
infer any components of $\mathcal{R}(\mathbf{C}_1^*)$ orthogonal to
$\mathcal{R}(\mathbf{Z})$ since they are annihilated by
$\mathbf{Z}^*$. Therefore, we apply Occam's razor and assume that
$\mathbf{C}_1^*=\mathbf{U}_1\hat{\mathbf{C}}_1^*$ for some
$\hat{\mathbf{C}}_1^*\in\mathbb{C}^{r_1\times n}$. By the same argument,
$\mathcal{R}(\mathbf{C}_2^*)$ cannot have any components orthogonal to
$\mathcal{R}(\mathbf{\Psi}_{\mathbf{Z}})$ since they are annihilated by
$\mathbf{\Psi}_{\mathbf{Z}}^*$ and have a positive contribution to the
regularization penalty term $\Tr (\mathbf{C}_2\mathbf{C}_2^*)$. Hence, we must
also have $\mathbf{C}_2^*=\mathbf{U}_2\hat{\mathbf{C}}_2^*$ for some
$\hat{\mathbf{C}}_2^*\in\mathbb{C}^{r_2\times n}$. Substituting these
relationships into \cref{eqn:multikernel_regression_loss_1} allows it to be
formulated as the standard least squares problem
\begin{equation}
  \begin{aligned}
    J &= \left\Vert \mathbf{X}^* -
      \mathbf{V}_1\mathbf{\Sigma}_1\hat{\mathbf{C}}_1^*
      - \mathbf{V}_2\mathbf{\Sigma}_2\hat{\mathbf{C}}_2^* \right\Vert_F^2 +
    \gamma \left\Vert \hat{\mathbf{C}}_2 \right\Vert_F^2 \\
    &= \left\Vert
      \begin{bmatrix}
        \mathbf{X}^* \\
        \mathbf{0}_{r_2\times n}
      \end{bmatrix} -
      \begin{bmatrix}
        \mathbf{V}_1\mathbf{\Sigma}_1 & \mathbf{V}_2\mathbf{\Sigma}_2 \\
        \mathbf{0}_{r_2\times r_1} & \sqrt{\gamma}\mathbf{I}_{r_2}
      \end{bmatrix}
      \begin{bmatrix} \hat{\mathbf{C}}_1^* \\
        \hat{\mathbf{C}}_2^*
      \end{bmatrix} \right\Vert_F^2
  \end{aligned}.
\end{equation}
The block-wise matrix clearly has full column rank $r_1+r_2$ for
$\gamma > 0$ and the normal equation for this least squares problem are found by
projecting onto its range. The solution,
\begin{equation}
  \begin{bmatrix} \hat{\mathbf{C}}_1^* \\ \hat{\mathbf{C}}_2^* \end{bmatrix} =
  \begin{bmatrix} \mathbf{\Sigma}_1 & \mathbf{0}_{r_1\times r_2} \\ \mathbf{0}_{r_2\times r_1} & \mathbf{\Sigma}_2 \end{bmatrix}^{-1}
  \begin{bmatrix} \mathbf{I}_{r_1} & \mathbf{V}_1^*\mathbf{V}_2 \\ \mathbf{V}_2^*\mathbf{V}_1 & \mathbf{I}_{r_2} + \gamma \mathbf{\Sigma}_2^{-2} \end{bmatrix}^{-1} \begin{bmatrix} \mathbf{V}_1^*\mathbf{X}^* \\ \mathbf{V}_2^*\mathbf{X}^* \end{bmatrix},
\end{equation}
corresponds to taking the left
pseudoinverse and simplifying the resulting expression. The matrices $\mathbf{V}_{1,2}$ and $\mathbf{\Sigma}_{1,2}$ are found by solving Hermitian eigenvalue problems using the (kernel)
matrices of inner products
$\mathbf{Z}^*\mathbf{Z} = \mathbf{V}_1\mathbf{\Sigma}_1^2\mathbf{V}_1^*$ and
$\mathbf{\Psi}_{\mathbf{Z}}^*\mathbf{\Psi}_{\mathbf{Z}} =
\mathbf{V}_2\mathbf{\Sigma}_2^2\mathbf{V}_2^*$. At a new point $\mathbf{z}$, the
approximate reconstruction using the partially linear kernel regression model is
\begin{equation}
  \label{eqn:partially_linear_regressor}
  \mathbf{x} \approx \left(\hat{\mathbf{C}}_1\mathbf{\Sigma}_1^{-1}\mathbf{V}_1^*\mathbf{Z}^*\right)\mathbf{z} + \hat{\mathbf{C}}_2\mathbf{\Sigma}_2^{-1}\mathbf{V}_2^*\left(\mathbf{\Psi}_{\mathbf{Z}}^*\mathbf{\Psi}(\mathbf{z})\right).
\end{equation}
Recall that the kernel matrices are
\begin{equation}
  \mathbf{\Psi}_{\mathbf{Z}}^*\mathbf{\Psi}_{\mathbf{Z}} =
  \begin{bmatrix}
    k(\mathbf{z}_1,\mathbf{z}_1) & \cdots & k(\mathbf{z}_1,\mathbf{z}_M) \\
    \vdots & \ddots & \vdots \\
    k(\mathbf{z}_M,\mathbf{z}_1) & \cdots & k(\mathbf{z}_M,\mathbf{z}_M)
  \end{bmatrix}, \qquad
  \mathbf{\Psi}_{\mathbf{Z}}^*\mathbf{\Psi}(\mathbf{z}) =
  \begin{bmatrix}
    k(\mathbf{z}_1,\mathbf{z}) \\
    \vdots  \\
    k(\mathbf{z}_M,\mathbf{z})
  \end{bmatrix}
\end{equation}
for a chosen continuous nonnegative definite mercer kernel function
$k:\mathbb{C}^d\times\mathbb{C}^d\rightarrow\mathbb{C}$ inducing the feature map
$\mathbf{\Psi}$.

The main drawback associated with the kernel method used for encoding
eigenfunctions or reconstructing the state is the number of kernel
evaluations. Even if the dimension $d$ of the reduced order model is small, the
kernel-based inner product of each new example must still be computed with all
of the training data in order to encode it and then again to decode it. When the
training data sets grow large, this leads to a high cost in making predictions
on new data points. An important avenue of future research is to prune the
training examples to only a small number of maximally informative ``support
vectors'' for taking inner products. Some possible approaches are discussed in
\cite{Tao2009,Peng2012}.

\section{Numerical examples}
\subsection{Duffing equation}
In our first numerical example, we will consider the unforced Duffing equation
in a parameter regime exhibiting two stable spirals. We take this example
directly from \cite{Williams2015} where the Koopman eigenfunctions are used to
separate and parameterize the basins of attraction for the fixed points. The
unforced Duffing equation is given by
\begin{equation}
  \label{eqn:Duffing}
  \ddot{x} = -\delta \dot{x} - x(\beta + \alpha x^2),
\end{equation}
where the parameters $\delta=0.5$, $\beta=-1$, and $\alpha=1$ are
chosen. The equation exhibits stable equilibria at $x=\pm 1$ with eigenvalues
$\lambda_{1,2}=\dfrac{1}{4}\left(-1 \pm \sqrt{31}\imath \right)$ associated with
the linearizations at these points. One can show that these (continuous-time)
eigenvalues also correspond to Koopman eigenfunctions whose magnitude and
complex argument act like action-angle variables parameterizing the entire
basins. A non-trivial Koopman eigenfunction with eigenvalue $\lambda_0 = 0$
takes different constant values in each basin, acting like an indicator function
to distinguish them.

We will see whether the LRAN and the reduced KDMD model can learn these
eigenfunctions from data and use them to predict the dynamics of the unforced
Duffing equation as well as to determine which basin of attraction a given point
belongs. The training data are generated by simulating the unforced Duffing
equation from uniform random initial conditions
$\left(x(0), \dot{x}(0)\right) \in [-2,2]\times [-2,2]$. From each trajectory
$11$ samples are recorded $\Delta t=0.25$ apart. The training data for LRAN
models consists of $M=10^4$ such trajectories. Since the KDMD method requires us
to evaluate the kernel function between a given example and each training point,
we limit the number of training data points to $10^3$ randomly selected snapshot
pairs from the original set. It is worth mentioning that the LRAN model handles
large data sets more efficiently than KDMD since the significant cost goes into
training the model which is then inexpensive and fast to evaluate on new
examples.

Since three of the Koopman eigenvalues are known ahead of time we train an LRAN
model where the transition matrix $\mathbf{K}$ is fixed to have discrete time
eigenvalues $\mu_k = \exp (\lambda_k \Delta t)$. We refer to this as the
``constrained LRAN'' and compare its performance to a ``free LRAN'' model where
$\mathbf{K}$ is learned and a $5$th order balanced truncation using KDMD called
``KDMD ROM''. The hyperparameters of each model are reported in
\cref{app:Hyperparameters}.

\begin{figure}[htbp]
  \centering
  \subfloat{
    \includegraphics[width=0.333\textwidth]{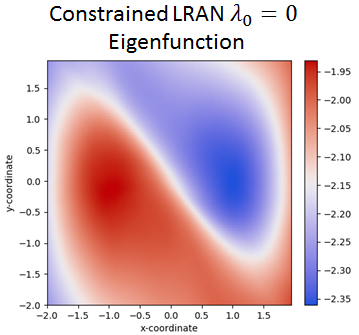}
  }
  \hfill
  \subfloat{
    \includegraphics[width=0.62\textwidth]{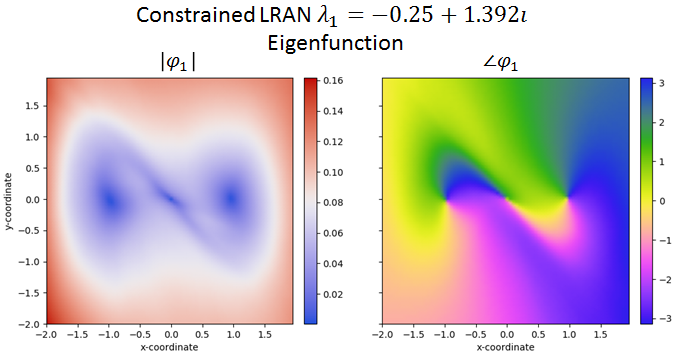}
  }
  \caption{Unforced Duffing eigenfunctions learned using constrained LRAN}
  \label{fig:Duffing_ConstrainedLRAN_Eigenfunctions}
\end{figure}

\begin{figure}[htbp]
  \centering
  \subfloat{
    \includegraphics[width=0.33\textwidth]{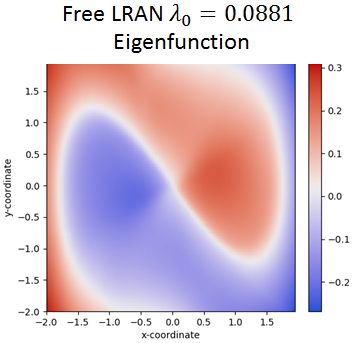}
  }
  \hfill
  \subfloat{
    \includegraphics[width=0.62\textwidth]{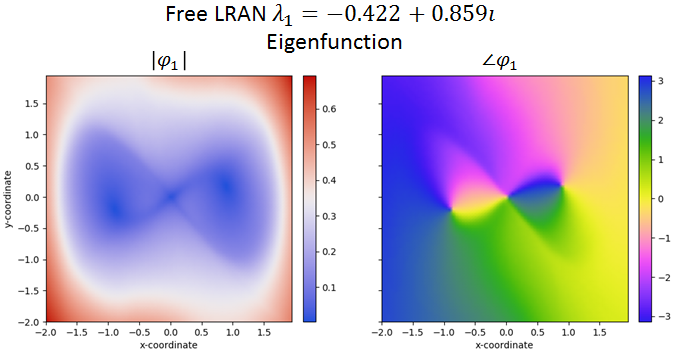}
  }
  \caption{Unforced Duffing eigenfunctions learned using free LRAN}
  \label{fig:Duffing_FreeLRAN_Eigenfunctions}
\end{figure}

\begin{figure}[htbp]
  \centering
  \subfloat{
    \includegraphics[width=0.325\textwidth]{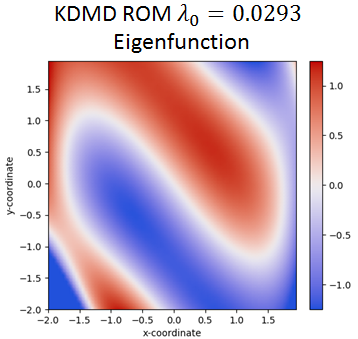}
  }
  \hfill
  \subfloat{
    \includegraphics[width=0.64\textwidth]{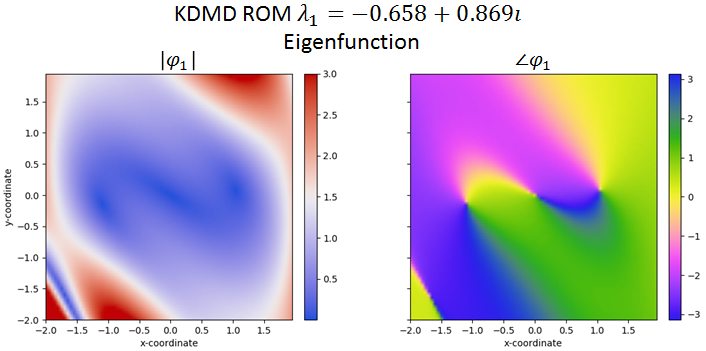}
  }
  \caption{Unforced Duffing eigenfunctions found using KDMD balanced ROM}
  \label{fig:Duffing_KDMDROM_Eigenfunctions}
\end{figure}

The learned eigenfunctions for each model are plotted in
\cref{fig:Duffing_ConstrainedLRAN_Eigenfunctions,fig:Duffing_FreeLRAN_Eigenfunctions,fig:Duffing_KDMDROM_Eigenfunctions}. The
corresponding eigenvalues learned or fixed in the model are also reported. The complex eigenfunctions are plotted in terms of their magnitude and phase. In
each case, the eigenfunction associated with the continuous-time eigenvalue
$\lambda_0$ closest to zero appears to partition the phase space into basins of
attraction for each fixed point as one would expect. In order to test this
hypothesis, we use the median eigenfunction value for each model as a threshold
to classify test data points between the basins. The eigenfunction learned by
the constrained LRAN was used to correctly classify $0.9274$ of the testing data
points. The free LRAN eigenfunction and the KDMD balanced reduced order model
eigenfunction correctly classified $0.9488$ and $0.9650$ of the testing data
respectively. $M_{\mbox{test}}=11*10^4$ test data points were used to evaluate
the LRAN models, though this number was reduced to a randomly selected
$M_{\mbox{test}}=1000$ to test the KDMD model due to the exceedingly high
computational cost of the kernel evaluations.

\begin{figure}[htbp]
  \centering
  \includegraphics[width=0.5\textwidth]{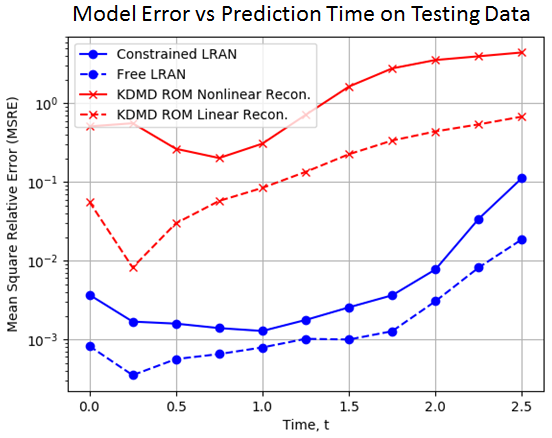}
  \caption{Unforced Duffing testing data mean square relative prediction errors
    for each model plotted against the prediction time}
  \label{fig:Duffing_ModelErrors}
\end{figure}

The other eigenfunction learned in each case parameterizes the basins of
attraction and therefore is used to account for the dynamics in each basin. Each
model appears to have learned a similar action-angle parameterization regardless
of whether the eigenvalues were specified ahead of time. However, the
constrained LRAN shows the best agreement with the true fixed point locations at
$x=\pm 1$ where $\left|\varphi_1\right| \rightarrow 0$. The mean square relative
prediction error was evaluated for each model by making predictions on the
testing data set at various times in the future. The results plotted in
\cref{fig:Duffing_ModelErrors} show that the free LRAN has by far the lowest
prediction error likely due to the lack of constraints on the functions it could
learn. It is surprising however, that nonlinear reconstruction hurt the
performance of the KDMD reduced order model. This illustrates a potential
difficulty with this method since the nonlinear part of the reconstruction is
prone to over-fit without sufficient regularization.

\subsection{Cylinder wake}
\begin{figure}[htbp]
  \centering
  \subfloat[]{
    \includegraphics[width=0.48\textwidth]{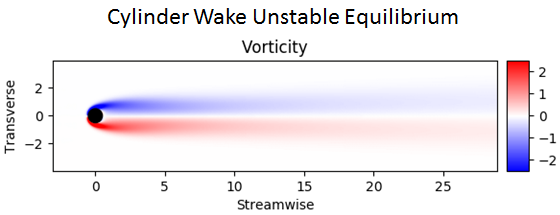}
    \label{fig:unstableEquilibrium}
  }
  \hfill
  \subfloat[]{
    \includegraphics[width=0.48\textwidth]{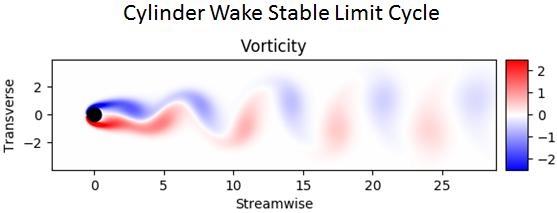}
    \label{fig:stableLimitCycle}
  }
  \label{fig:cylinderWake_exampleSnapshots}
  \caption{Example cylinder wake flow snapshots at the unstable equilibrium and on the stable limit cycle}
\end{figure}

The next example we consider is the formation of a K\'{a}rm\'{a}n vortex sheet
downstream of a cylinder in a fluid flow. This problem was chosen since the data
has low intrinsic dimensionality due to the simple flow structure but is
embedded in high-dimensional snapshots. We are interested in whether the
proposed techniques can be used to discover very low dimensional models that
accurately predict the dynamics over many time steps. We consider the growth of
instabilities near an unstable base flow shown in \cref{fig:unstableEquilibrium} at Reynold number $Re=60$
all the way until a stable limit cycle shown in \cref{fig:stableLimitCycle} is reached. The models will have to learn
to make predictions over a range of unsteady flow conditions from the unstable
equilibrium to the stable limit cycle.

The raw data consisted of $2000$ simulated snapshots of the velocity field taken at time intervals
$0.2 D/U_{\infty}$, where $D$ is the cylinder diameter and $U_{\infty}$ is the
free-stream velocity. These data were split into $M_{\mbox{train}}=1000$
training, $M_{\mbox{eval}}=500$ evaluation, and $M_{\mbox{test}}=500$ testing
data points. Odd numbered points $\Delta t = 0.4 D/U_{\infty}$ apart were used
for training. The remaining $1000$ points were divided again into even and odd
numbered points $2\Delta t = 0.8 D/U_{\infty}$ apart for evaluation and
testing. This enabled training, evaluation, and testing on long data sequences
while retaining coverage over the complete trajectory. The continuous-time eigenvalues are found from the discrete-time eigenvalues according to $\lambda = \log(\mu)/\Delta t = \log(\mu)U_{\infty}/(0.4D)$.

\begin{figure}[htbp]
  \centering
  \subfloat[]{
    \includegraphics[width=0.48\textwidth]{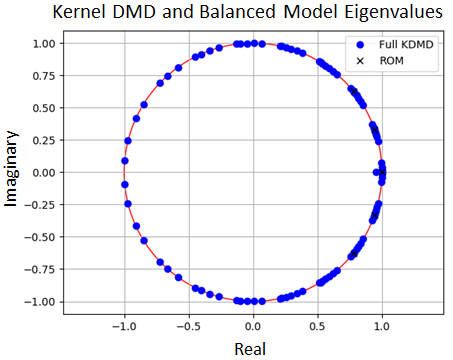}
    \label{fig:cylinderWake_KDMDROM_Eigenvalues}
  }
  \hfill
  \subfloat[]{
    \includegraphics[width=0.48\textwidth]{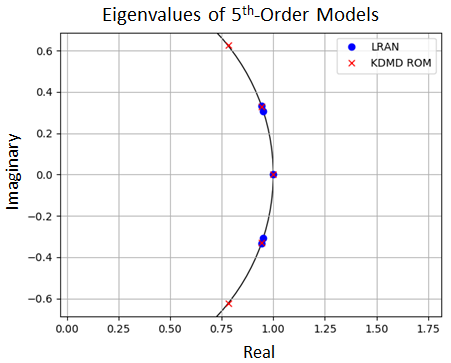}
    \label{fig:cylinderWake_Models_Eigenvalues}
  }
  \label{fig:cylinderWake_KoopmanEigenvalues}
  \caption{Discrete-time Koopman eigenvalues approximated by the KDMD ROM and the LRAN}
\end{figure}

The raw data was projected onto its $200$ most energetic POD modes which
captured essentially all of the energy in order to reduce the cost of storage
and training. $400$-dimensional time delay embedded snapshots were formed from
the state at time $t$ and $t+\Delta t$. A $5$th-order LRAN model and the
5th-order KDMD reduced order model were trained using the hyperparameters in
\cref{tab:CylinderWake_LRAN_Hyperparameters,tab:CylinderWake_KDMDROM_Hyperparameters}. In
\cref{fig:cylinderWake_KDMDROM_Eigenvalues}, many of the discrete-time eigenvalues
given by the over-specified KDMD model have approximately neutral stability with some being slightly unstable. However, the finite horizon
formulation for balanced truncation allows us to learn the most dynamically
salient eigenfunctions over a given length of time, in this case
$\mathcal{T}=20$ steps or $8.0 D/U_{\infty}$. We see in
\cref{fig:cylinderWake_Models_Eigenvalues} that three of the eigenvalues learned
by the two models are in close agreement and all are approximately neutrally
stable.

\begin{figure}[htbp]
  \centering
  \subfloat{
    \includegraphics[width=0.48\textwidth]{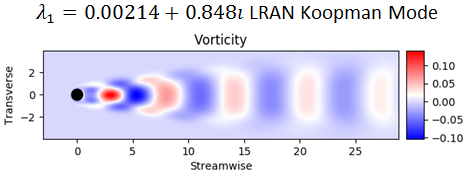}
  }
  \hfill
  \subfloat{
    \includegraphics[width=0.48\textwidth]{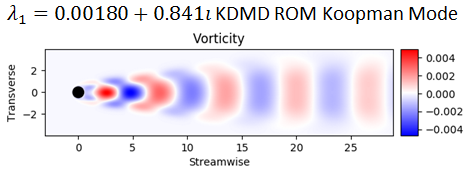}
  }
  \label{fig:cylinderWake_KoopmanModes_1}
  \caption{LRAN and KDMD ROM Koopman modes associated with
    $\lambda_1 \approx 0.002 + 0.845\imath$}
\end{figure}

\begin{figure}[htbp]
  \centering
  \subfloat{
    \includegraphics[width=0.48\textwidth]{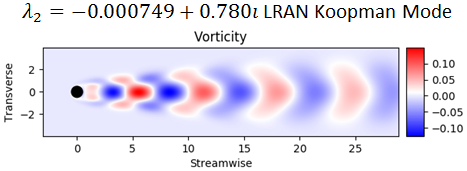}
  }
  \hfill
  \subfloat{
    \includegraphics[width=0.48\textwidth]{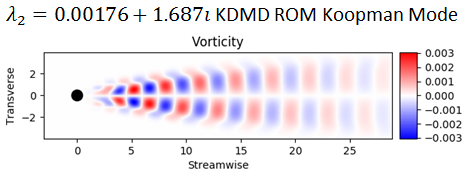}
  }
  \label{fig:cylinderWake_KoopmanModes_2}
  \caption{LRAN and KDMD ROM Koopman modes associated with $\lambda_2$ which
    differs greatly between the models}
\end{figure}

A side-by-side comparison of the Koopman modes gives some insight into the flow structures whose dynamics the Koopman eigenvalues describe.
We notice right away that the Koopman modes in \cref{fig:cylinderWake_KoopmanModes_1} corresponding to continuous-time eigenvalue $\lambda_1$ are very similar for both models and indicate the pattern of vortex shedding downstream. 
This makes sense since a single frequency and mode will account for most of the amplitude as the limit cycle is approached. 
Evidently both models discover these limiting periodic dynamics. 
For the KDMD ROM, $\lambda_2$ is almost exactly the higher harmonic $2*\lambda_1$. 
The corresponding Koopman mode in \cref{fig:cylinderWake_KoopmanModes_2} also reflects smaller flow structures which oscillate at twice the frequency of $\lambda_1$. 
Interestingly, the LRAN does not learn the same second eigenvalue as the KDMD ROM. 
The LRAN continuous-time eigenvalue $\lambda_2$ is very close to $\lambda_1$ which suggest that these frequencies might team up to produce the low-frequency $\lambda_1 - \lambda_2$. 
The second LRAN Koopman mode in \cref{fig:cylinderWake_KoopmanModes_2} also bears qualitative resemblance to the first Koopman mode in \cref{fig:cylinderWake_KoopmanModes_1}, but with a slightly narrower pattern in the y-direction. 
The LRAN may be using the information at these frequencies to capture some of the slower transition process from the unstable fixed point to the limit cycle. 
The Koopman modes corresponding to $\lambda_0=0$ are also qualitatively different indicating that the LRAN and KDMD ROM are extracting different constant features from the data.
We must be careful in our interpretation, however, since the LRAN's koopman modes are only a least squares approximations to the nonlinear reconstruction process performed by the decoder.

\begin{figure}[htbp]
  \centering
  \subfloat{
    \includegraphics[width=0.48\textwidth]{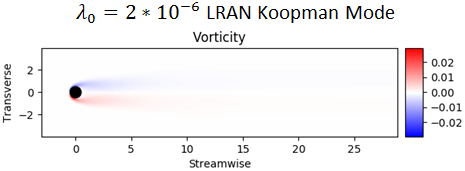}
  }
  \hfill
  \subfloat{
    \includegraphics[width=0.48\textwidth]{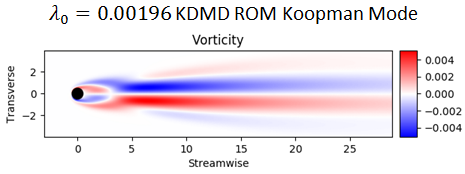}
  }
  \label{fig:cylinderWake_KoopmanModes_0}
  \caption{LRAN and KDMD ROM Koopman modes associated with $\lambda_0 \approx 0$}
\end{figure}

\begin{figure}[htbp]
  \centering
  \includegraphics[width=0.5\textwidth]{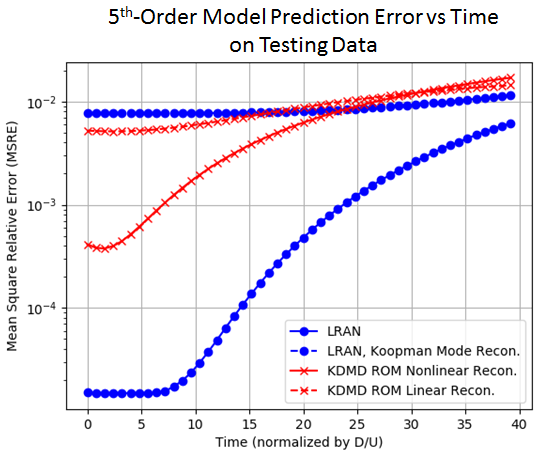}
  \caption{Cylinder wake testing data mean square relative prediction errors for
    each model plotted against the prediction time}
  \label{fig:CylinderWake_ModelErrors}
\end{figure}

Plotting the model prediction error \cref{fig:CylinderWake_ModelErrors} shows
that the linear reconstructions using both models have comparable performance
with errors growing slowly over time. Therefore, the choice of the second
Koopman mode does not seem to play a large role in the reconstruction
process. However, when the nonlinear decoder is used to reconstruct the LRAN
predictions, the mean relative error is roughly an order of magnitude smaller
than the nonlinearly reconstructed KDMD ROM over many time steps. The LRAN has
evidently learned an advantageous nonlinear transformation for reconstructing
the data using the features evolving according to $\lambda_2$. The second
Koopman mode reflects a linear approximation of this nonlinear transformation in
the least squares sense.

Another remark is that nonlinear reconstruction using the KDMD ROM did
significantly improve the accuracy in this example. This indicates that many of
the complex variations in the data are really enslaved to a small number of
modes. This makes sense since the dynamics are periodic on the limit cycle. 
Finally, it is worth mentioning that the prediction accuracy was
achieved on average over all portions of the trajectory from the unstable
equilibrium to the limit cycle. Both models therefore have demonstrated
predictive accuracy and validity over a wide range of qualitatively different
flow conditions. The nonlinearly reconstructed LRAN achieves a constant low
prediction error over the entire time interval used for training
$\mathcal{T}\Delta t = 8.0 D/U_{\infty}$. The error only begins to grow outside
the interval used for training. The high prediction accuracy could likely be
extended by training on longer data sequences.

\subsection{Kuramoto-Sivashinsky equation}
We now move on to test our new techniques on a very challenging problem --- the
Kuramoto-Sivashinsky equation in a parameter regime just beyond the onset of
chaos. Since any chaotic dynamical system is mixing, it only has trivial Koopman
eigenfunctions on its attractor(s). We therefore cannot expect our model to 
accurately reflect the dynamics of the real system. Rather, we aim to make
predictions using very low-dimensional models that are accurate over short
times and plausible over longer times.

The data was generated by performing direct numerical simulations of the
Kuramoto-Sivashinsky equation, 
\begin{equation}
\label{eqn:KuramotoSivashinksy}
u_t + u_{xx} + u_{xxxx} + u u_x = 0, \qquad x\in [0,L],
\end{equation}
using a semi-implicit Fourier pseudo-spectral method. The length $L=8\pi$ was chosen
where the equation first begins to exhibit chaotic dynamics
\cite{Holmes2012}. $128$ Fourier modes were used to resolve all of the
dissipative scales. Each data set: training, evaluation, and test, consisted of
$20$ simulations from different initial conditions each with $500$ recorded
states spaced by $\Delta t = 1.0$. Snapshots consisted of time delay embedded
states at $t$ and $t+\Delta t$. The initial conditions were formed by suppling
Gaussian random perturbations to the coefficients on the $3$ linearly unstable
Fourier modes $0 < 2\pi k/L < 1 \Longrightarrow k=1,2,3$.

\begin{figure}[htbp]
  \centering
  \subfloat[]{
    \includegraphics[width=0.48\textwidth]{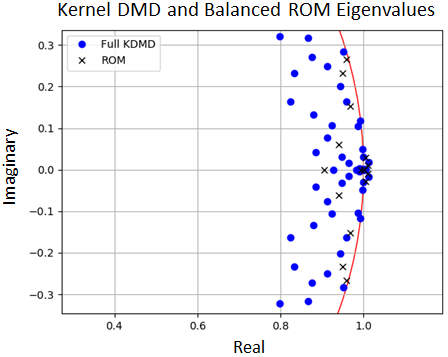}
    \label{fig:KuramotoSivashinsky_KDMDROM_Eigenvalues}
  }
  \hfill
  \subfloat[]{
    \includegraphics[width=0.48\textwidth]{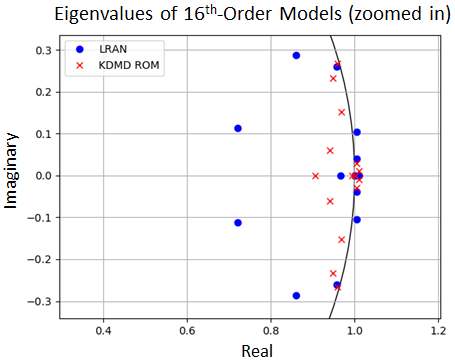}
    \label{fig:KuramotoSivashinsky_Models_Eigenvalues}
  }
  \label{fig:KuramotoSivashinksy_KoopmanEigenvalues}
  \caption{Discrete-time Koopman eigenvalues for the Kuramoto-Sivashinksy
    equation approximated by the KDMD ROM and the LRAN}
\end{figure}

An LRAN as well as a KDMD balanced ROM were trained to make predictions over a
time horizon $\mathcal{T}=5$ steps using only $d=16$ dimensional models. Model
parameters are given in
\cref{tab:KuramotoSivashinksy_LRAN_Hyperparameters,tab:KuramotoSivashinksy_KDMDROM_Hyperparameters}. The
learned approximate Koopman eigenvalues are plotted in
\cref{fig:KuramotoSivashinksy_KoopmanEigenvalues}. We notice that there are some
slightly unstable eigenvalues, which makes sense since there are certainly
unstable modes including the three linearly unstable Fourier
modes. Additionally, \cref{fig:KuramotoSivashinsky_Models_Eigenvalues} shows
that some of the eigenvalues with large magnitude learned by the LRAN and the KDMD ROM
are in near agreement.

\begin{figure}[htbp]
  \centering
  \includegraphics[width=0.5\textwidth]{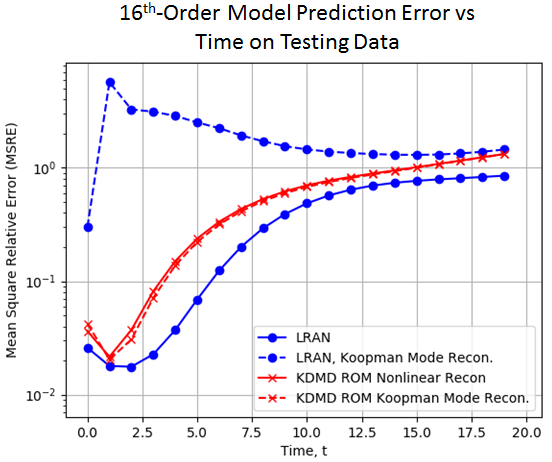}
  \caption{Kuramoto-Sivashinsky testing data mean square relative prediction errors for each model plotted against the prediction time}
  \label{fig:KuramotoSivashinsky_ModelErrors}
\end{figure}

The plot of mean square relative prediction error on the testing data set
\cref{fig:KuramotoSivashinsky_ModelErrors} indicates that our addition of
nonlinear reconstruction from the low dimensional KDMD ROM state does not change
the accuracy of the reconstruction. The performance of the KDMD ROM and the LRAN
are comparable with the LRAN showing a modest reduction in error over all
prediction times. It is interesting to note that the LRAN does not produce
accurate reconstructions using the regression-based Koopman modes. In this
example, the LRAN's nonlinear decoder is essential for the reconstruction
process. Evidently, the dictionary functions learned by the encoder require
nonlinearity to reconstruct the state. Again, both models are most accurate over
the specified time horizon $\mathcal{T}=5$ used for training.

\begin{figure}[htbp]
  \centering
  \subfloat[]{
    \includegraphics[width=0.485\textwidth]{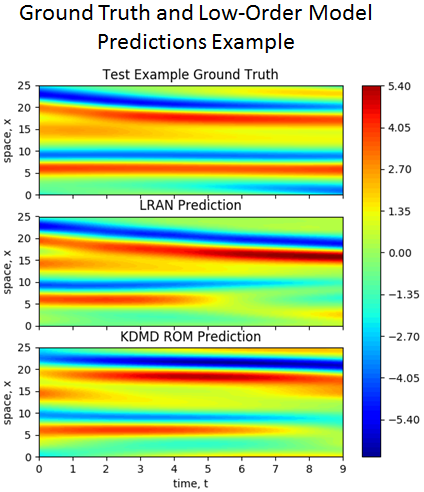}
    \label{fig:KuramotoSivashinsky_ExamplePred_1}
  }
  \hfill
  \subfloat[]{
    \includegraphics[width=0.475\textwidth]{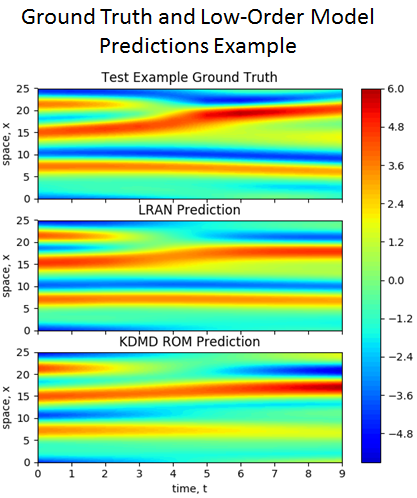}
    \label{fig:KuramotoSivashinsky_ExamplePred_2}
  }
  \label{fig:KuramotoSivashinsky_ExamplePreds}
  \caption{LRAN and KDMD ROM model predictions on Kuramoto-Sivashinsky test data examples}
\end{figure}

Plotting some examples in \cref{fig:KuramotoSivashinsky_ExamplePreds} of ground
truth and predicted test data sequences illustrates the behavior of the
models. These examples show that both the LRAN and the KDMD ROM make
quantitatively accurate short term predictions. While the predictions after
$t \approx 5$ lose their accuracy as one would expect when trying to make linear
approximations of chaotic dynamics, they remain qualitatively plausible. The
LRAN model in particular is able to model and predict grouping and merging
events between traveling waves in the solution. For example in
\cref{fig:KuramotoSivashinsky_ExamplePred_1} the LRAN successfully predicts the
merging of two wave crests (in red) taking place between $t=2$ and $t=5$. The
LRAN also predicts the meeting of a peak and trough in
\cref{fig:KuramotoSivashinsky_ExamplePred_2} at $t=5$. These results are
encouraging considering the substantial reduction in dimensionality from a time
delay embedded state of dimension $256$ to a $16$-dimensional encoded state
having linear time evolution.

\section{Conclusions}
We have illustrated some fundamental
challenges with EDMD, in particular highlighting the trade-off between rich dictionaries and
over-fitting. The use of adaptive, low-dimensional dictionaries avoids the over-fitting problem while retaining enough capacity to represent
Koopman eigenfunctions of complicated systems. This motivates the use of neural
networks to learn sets of dictionary observables that are well-adapted to the
given problems. By relaxing the constraint that the models must produce linear
reconstructions of the state via the Koopman modes, we introduce a decoder
neural network enabling the formation of very low-order models utilizing richer
observables. Finally, by combining the neural network architecture that is
essentially an autoencoder with linear recurrence, the LRAN learns features
that are dynamically important rather than just energetic (i.e., large in norm).

Discovering a small set of dynamically important features or states is also the
idea behind balanced model reduction. This led us to investigate the
identification of low-dimensional models by balanced truncation of
over-specified EDMD models, and in particular, KDMD models. Nonlinear reconstruction
using a partially linear multi-kernel method was investigated for improving the
reconstruction accuracy of the KDMD ROMs from very low-dimensional spaces. Our
examples show that in some cases like the cylinder wake example, it can greatly
improve the accuracy. We think this is because the data is intrinsically
low-dimensional, but curves in such a way as to extend in many dimensions of the
embedding space. The limiting case of the cylinder flow is an extreme example
where the data becomes one-dimensional on the limit cycle. In some other cases,
however, nonlinear reconstruction does not help, is sensitive to parameter
choices, or decreases the accuracy due to over-fitting.

Our numerical examples indicate that unfolding the linear recurrence for many
steps can improve the accuracy of LRAN predictions especially within the time
horizon used during training. This is observed in the error versus prediction
time plots in our examples: the error remains low and relatively flat for
predictions made inside the training time horizon $\mathcal{T}$. The error then
grows for predictions exceeding this length of time. However, for more
complicated systems like the Kuramoto-Sivashinsky equation, one cannot unfold
the network for too many steps before additional dimensions must be added to
retain accuracy of the linear model approximation over time. These observations
are also approximately true of the finite-horizon BPOD formulation used to
create approximate balanced truncations of KDMD models. One additional
consideration in forming balanced reduced-order models from finite-horizon
impulse responses of over-specified systems is the problem of spurious
eigenvalues whose associated modes only become significant for approximations as
$t \rightarrow \infty$. The use of carefully chosen finite time horizons allows
us to pick features which are the most relevant (observable and excitable) over
any time span of interest.

The main drawback of the balanced reduced-order KDMD models becomes evident when
making predictions on new evaluation and testing examples. While the LRAN has a
high ``up-front'' cost to train --- typically requiring hundreds of thousands of
iterations --- the cost of evaluating a new example is almost negligible and so
very many predictions can be made quickly. On the other hand, every new example
whose dynamics we want to predict using the KDMD reduced order model must still
have its inner product evaluated with every training data point. Pruning methods
like those developed for multi-output least squares support vector machines
(LS-SVMs) will be needed to reduce the number of kernel evaluations before KDMD
reduced order models can be considered practical for the purpose of predicting
dynamics. The same will need to be done for kernel-based nonlinear
reconstruction methods.

We conclude by discussing some exciting directions for future work on LRANs. The
story certainly does not end with using them to learn linear encoded state
dynamics and approximations of the Koopman eigenfunctions. The idea of
establishing a possibly complicated transformation into and out of a space where
the dynamics are simple is an underlying theme of this work. Human understanding
seems to at least partly reside in finding isomorphism. If we can establish that
a seemingly complicated system is topologically conjugate to a simple system, it
doesn't really matter what the transformation is as long as we can compute it. In
this vein, the LRAN architecture is perfectly suited to learning the complicated
transformations into and out of a space where the dynamics equations have a
simple normal form. For example, this could be accomplished by learning
coefficients on homogeneous polynomials of varying degree in addition to a
matrix for updating the encoded state. One could also include learned parameter
dependencies in order to study bifurcation behavior. Furthermore, any kind of
symmetry could be imposed through either the normal form itself or through the
neural network's topology. For example, convolutional neural networks can be
used in cases where the system is statistically stationary in a spatial
coordinate.

Introduction of control terms to the dynamics of the encoded state is another
interesting direction for inquiry. In some cases it might be possible to
introduce a second autoencoder to perform a state-dependent encoding and
decoding of the control inputs at each time step. Depending on the form that the
inputs take in the evolution of the encoded state, it may be possible to apply a
range of techniques from modern state-space control to nonlinear optimal control
or even reinforcement learning-based control with policy and value functions
parameterized by neural networks.

Another natural question is how the LRAN framework can be adapted to handle
stochastic dynamics and uncertainty. Recent work in the development of
structured inference networks for nonlinear stochastic models
\cite{Krishnan2017} may offer a promising approach. The low-dimensional dynamics
and reconstruction process could be generalized to nonlinear stochastic
processes for generating full state trajectories. Since the inference problem
for such nonlinear systems is intractable, the encoder becomes a
(bi-directional) recurrent neural network for performing approximate inference
on the latent state given our data in analogy with Kalman smoothing. In this
manner, many plausible outputs can be generated to estimate the distribution
over state trajectories in addition to an inference distribution for quantifying
uncertainty about the low-dimensional latent state.

Furthermore, with the above formulation of generative autoencoder networks for
dynamics, it might be possible to employ adversarial training
\cite{goodfellow2014} in a similar manner to the adversarial autoencoder
\cite{Makhzani2015}. Training the generative network used for reconstruction
against a discriminator network will encourage the generator to produce more
plausible details like turbulent eddies in fluid flows which are not easily
distinguished from the real thing.

\appendix
\section{Hyperparameters used to train models}
\label{app:Hyperparameters}
The same hyperparameters in \cref{tab:UnforcedDuffing_LRAN_Hyperparameters} were
used to train the constrained and free LRANs in the unforced Duffing equation
example.
\begin{table}[htbp]
  \caption{Constrained LRAN hyperparameters for unforced Duffing example}
  \label{tab:UnforcedDuffing_LRAN_Hyperparameters}
  \centering
  \begin{tabular}{|c|c|} \hline
   \bf Parameter & \bf Value(s) \\ \hline
    Time-delays embedded in a snapshot & $1$ \\
    Encoder layer widths (left to right) & $2$, $32$, $32$, $16$, $16$, $8$, $3$ \\
    Decoder layer widths (left to right) & $3$, $8$, $16$, $16$, $32$, $32$, $2$ \\
    Snapshot sequence length, $\mathcal{T}$ & $10$ \\
    Weight decay rate, $\delta$ & $0.8$ \\
    Relative weight on encoded state, $\beta$ & $1.0$ \\
    Minibatch size & $50$ examples \\
    Initial learning rate & $10^{-3}$ \\
    Geometric learning rate decay factor & $0.01$ per $4*10^5$ steps \\
    Number of training steps & $4*10^5$ \\
    \hline
  \end{tabular}
\end{table}

\Cref{tab:UnforcedDuffing_KDMDROM_Hyperparameters} summarizes the
hyperparameters used to train the KDMD Reduced Order Model for the unforced
Duffing example.

\begin{table}[htbp]
  \caption{KDMD ROM hyperparameters for unforced Duffing example}
  \label{tab:UnforcedDuffing_KDMDROM_Hyperparameters}
  \centering
  \begin{tabular}{|c|c|} \hline
   \bf Parameter & \bf Value(s) \\ \hline
    Time-delays embedded in a snapshot & $1$ \\
    EDMD Dictionary kernel function & Gaussian RBF, $\sigma = 10.0$ \\
    KDMD SVD rank, $r$ & $27$ \\
    BPOD time horizon, $\mathcal{T}$ & $10$ \\
    BPOD output projection rank & $2$ (no projection) \\
    Balanced model order, $d$ & $3$ \\
    Nonlinear reconstruction kernel function & Gaussian RBF, $\sigma = 10.0$ \\
    Multi-kernel linear part truncation rank, $r_1$ & $3$ \\
    Multi-kernel nonlinear part truncation rank, $r_2$ & $8$ \\
    Multi-kernel regularization constant, $\gamma$ & $10^{-4}$ \\
    \hline
  \end{tabular}
\end{table}

The hyperparameters used to train the LRAN model on the cylinder wake data are
given in \cref{tab:CylinderWake_LRAN_Hyperparameters}.
\begin{table}[htbp]
  \caption{LRAN hyperparameters for cylinder wake example}
  \label{tab:CylinderWake_LRAN_Hyperparameters}
  \centering
  \begin{tabular}{|c|c|} \hline
   \bf Parameter & \bf Value(s) \\ \hline
    Time-delays embedded in a snapshot & $2$ \\
    Encoder layer widths (left to right) & $400$, $100$, $50$, $20$, $10$, $5$ \\
    Decoder layer widths (left to right) & $5$, $10$, $20$, $50$, $100$, $400$ \\
    Snapshot sequence length, $\mathcal{T}$ & $20$ \\
    Weight decay rate, $\delta$ & $0.95$ \\
    Relative weight on encoded state, $\beta$ & $1.0$ \\
    Minibatch size & $50$ examples \\
    Initial learning rate & $10^{-3}$ \\
    Geometric learning rate decay factor & $0.01$ per $2*10^5$ steps \\
    Number of training steps & $2*10^5$ \\
    \hline
  \end{tabular}
\end{table}

\Cref{tab:CylinderWake_KDMDROM_Hyperparameters} summarizes the hyperparameters
used to train the KDMD Reduced Order Model on the cylinder wake data.

\begin{table}[htbp]
  \caption{KDMD ROM hyperparameters for cylinder wake example}
  \label{tab:CylinderWake_KDMDROM_Hyperparameters}
  \centering
  \begin{tabular}{|c|c|} \hline
   \bf Parameter & \bf Value(s) \\ \hline
    Time-delays embedded in a snapshot & $2$ \\
    EDMD Dictionary kernel function & Gaussian RBF, $\sigma = 10.0$ \\
    KDMD SVD rank, $r$ & $100$ \\
    BPOD time horizon, $\mathcal{T}$ & $20$ \\
    BPOD output projection rank & $100$ \\
    Balanced model order, $d$ & $5$ \\
    Nonlinear reconstruction kernel function & Gaussian RBF, $\sigma = 10.0$ \\
    Multi-kernel linear part truncation rank, $r_1$ & $5$ \\
    Multi-kernel nonlinear part truncation rank, $r_2$ & $15$ \\
    Multi-kernel regularization constant, $\gamma$ & $10^{-8}$ \\
    \hline
  \end{tabular}
\end{table}

\Cref{tab:KuramotoSivashinksy_LRAN_Hyperparameters} lists the hyperparameters
used to train the LRAN model on the Kuramoto-Sivashinsky equation example.
\begin{table}[htbp]
  \caption{LRAN hyperparameters for Kuramoto-Sivashinsky example}
  \label{tab:KuramotoSivashinksy_LRAN_Hyperparameters}
  \centering
  \begin{tabular}{|c|c|} \hline
   \bf Parameter & \bf Value(s) \\ \hline
    Time-delays embedded in a snapshot & $2$ \\
    Encoder layer widths (left to right) & $256$, $32$, $32$, $16$, $16$ \\
    Decoder layer widths (left to right) & $16$, $16$, $32$, $32$, $256$ \\
    Snapshot sequence length, $\mathcal{T}$ & $5$ \\
    Weight decay rate, $\delta$ & $0.9$ \\
    Relative weight on encoded state, $\beta$ & $1.0$ \\
    Minibatch size & $50$ examples \\
    Initial learning rate & $10^{-3}$ \\
    Geometric learning rate decay factor & $0.1$ per $2*10^5$ steps \\
    Number of training steps & $4*10^5$ \\
    \hline
  \end{tabular}
\end{table}

\Cref{tab:KuramotoSivashinksy_KDMDROM_Hyperparameters} summarizes the
hyperparameters used to train the KDMD Reduced Order Model on the
Kuramoto-Sivashinsky equation example data.

\begin{table}[htbp]
  \caption{KDMD ROM hyperparameters for Kuramoto-Sivashinsky example}
  \label{tab:KuramotoSivashinksy_KDMDROM_Hyperparameters}
  \centering
  \begin{tabular}{|c|c|} \hline
   \bf Parameter & \bf Value(s) \\ \hline
   	Time-delays embedded in a snapshot & $2$ \\
   	EDMD Dictionary kernel function & Gaussian RBF, $\sigma = 10.0$ \\
    KDMD SVD rank, $r$ & $60$ \\
    BPOD time horizon, $\mathcal{T}$ & $5$ \\
    BPOD output projection rank & $60$ \\
    Balanced model order, $d$ & $16$ \\
    Nonlinear reconstruction kernel function & Gaussian RBF, $\sigma = 100.0$ \\
    Multi-kernel linear part truncation rank, $r_1$ & $16$ \\
    Multi-kernel nonlinear part truncation rank, $r_2$ & $60$ \\
    Multi-kernel regularization constant, $\gamma$ & $10^{-7}$ \\
    \hline
  \end{tabular}
\end{table}

\section*{Acknowledgments}
We would like to thank Scott Dawson for providing us with the data from his
cylinder wake simulations. We would also like to thank William Eggert for his invaluable help and collaboration on the initial iterations of the LRAN code.

\bibliographystyle{siamplain}
\bibliography{RowleyGroupReferences}
\end{document}